\documentclass[12pt,bezier]{article}

\usepackage{amssymb,amsmath,amsthm,latexsym,tikz,url,subfig,float,amsthm,mathrsfs,thmtools,thm-restate}
\usepackage[normalem]{ulem}
\newtheorem{theorem}{Theorem}[section]
\newtheorem{corollary}[theorem]{Corollary}
\newtheorem{definition}[theorem]{Definition}
\newtheorem{conjecture}[theorem]{Conjecture}
\newtheorem{lemma}[theorem]{Lemma}

\newtheorem{remark}[theorem]{Remark}
\newtheorem{problem}[theorem]{Problem}

\newcommand{\y}{{\bf y}}
\newcommand{\1}{{\bf 1}}
\newcommand{\x}{{\bf x}}

\newcommand{\dm}{{\rm diam}}
\newcommand{\G}{{\cal G}}
\newcommand{\g}{\Gamma}
\newcommand{\e}{\epsilon}
\newcommand{\vertex}{\node[vertex]}
\tikzstyle{vertex}=[circle, draw, inner sep=0pt, minimum size=3pt]

\usepackage{color}

\begin{document}

\title{ Minimum algebraic connectivity and maximum diameter: Aldous--Fill and Guiduli--Mohar conjectures}
\author{ Maryam Abdi$^{\,\rm a}$ \quad  Ebrahim Ghorbani$^{\,\rm b,c,}$\thanks{Corresponding author, {\tt ebrahim.ghorbani@uni-hamburg.de} }
\\[.3cm]
{\sl\normalsize $^{\rm a}$School of Mathematics, Institute for Research in Fundamental Sciences (IPM),}\\
{\sl\normalsize P. O. Box 19395-5746, Tehran, Iran }\\
{\sl\normalsize $^{\rm b}$Department of Mathematics, K. N. Toosi University of Technology,}\\
{\sl\normalsize P. O. Box 16765-3381, Tehran, Iran}\\
{\sl\normalsize $^{\rm c}$Department of Mathematics, University of Hamburg,}\\
{\sl\normalsize Bundesstra\ss e 55 (Geomatikum), 20146 Hamburg, Germany }}

\maketitle
%\footnotetext{{\em E-mail Addresses}: {\tt m.abdi@email.kntu.ac.ir} (M. Abdi), {\tt ebrahim.ghorbani@uni-hamburg.de} (E. Ghorbani) }

\begin{abstract}
Aldous and Fill (2002) conjectured that the maximum relaxation time for
the random walk on a connected regular graph with $n$ vertices is
$(1+o(1)) \frac{3n^{2}}{2\pi ^{2}}$. A conjecture by Guiduli and Mohar
(1996) predicts the structure of graphs whose algebraic connectivity
$\mu $ is the smallest among all connected graphs whose minimum degree
$\delta $ is a given $d$. We prove that this conjecture implies the Aldous--Fill
conjecture for odd $d$. We pose another conjecture on the structure of
$d$-regular graphs with minimum $\mu $, and show that this also implies
the Aldous--Fill conjecture for even $d$. In the literature, it has been
noted empirically that graphs with small $\mu $ tend to have a large diameter.
In this regard, Guiduli (1996) asked if the cubic graphs with maximum diameter
have algebraic connectivity smaller than all others. Motivated by these,
we investigate the interplay between the graphs with maximum diameter and
those with minimum algebraic connectivity. We show that the answer to Guiduli
problem in its general form, that is for $d$-regular graphs for every
$d\ge 3$ is negative. We aim to develop an asymptotic formulation of the
problem. It is proven that 
$d$-regular graphs for $d\ge 5$ as well as graphs with
$\delta =d$ for $d\ge 4$  with asymptotically maximum diameter, do not necessarily exhibit the asymptotically smallest $\mu$. We conjecture that $d$-regular graphs (or graphs
with $\delta =d$) that have asymptotically smallest $\mu $, should have
asymptotically maximum diameter. 
The above results rely heavily on our
understanding of the structure as well as optimal estimation of the algebraic
connectivity of nearly maximum-diameter graphs, from which the Aldous--Fill conjecture for this family of graphs also follows.
\vspace{4mm}

\noindent {\bf AMS Mathematics Subject Classification\,(2010):}   05C50, 60G50, 05C35
 
\end{abstract}

\section{Introduction}
All graphs we consider are simple, i.e. undirected graphs without loops or multiple edges.
Additionally, we assume that they are connected.
The {\em relaxation time} of the random walk on a graph $G$  is defined by $\tau=1/(1-\eta_2)$, where $\eta_2$ is the second largest eigenvalue of
the {\em transition matrix} of $G$, that is the matrix $\Delta^{-1} A$ in which $\Delta$ and $A$ are the diagonal matrix of vertex degrees and the adjacency matrix of $G$, respectively.
A central problem in the study of random walks is to determine the {\em mixing time}, a measure of how fast the
random walk converges to the stationary distribution. As seen  through  the literature \cite{aldous2002reversible,chung}, the  relaxation time is the primary term controlling mixing time. Therefore, relaxation time is directly associated with the rate of convergence  of the random walk.
Our main motivation in this work is the following conjecture on the maximum relaxation time of the random walk in regular graphs.

\begin{conjecture}[{Aldous and Fill \cite[p.~217]{aldous2002reversible}}]\label{conj:A-F} \rm
	Over all  regular graphs on $n$ vertices, $\max \tau =(1+o(1)) \frac{3n^2}{2\pi^2}$.
\end{conjecture}

 For a graph $G$,  $L(G)=\Delta-A$ is its  Laplacian matrix. The second smallest eigenvalue of $L(G)$ is called the  {\em algebraic connectivity}  of $G$ and it is denoted by  $\mu=\mu(G)$.
When  $G$ is regular, of degree $d$ say, then its transition matrix is $\frac1d A$ and its Laplacian is  $dI-A$.
It is then seen that the relaxation time of $G$ is equal to $d/\mu(G)$. Also as $G$ is regular, $\mu(G)$ is the same as its {\em spectral gap}, the difference between the two largest eigenvalues of the adjacency matrix of  $G$.
So within the family of $d$-regular graphs, maximizing the relaxation time is equivalent to minimizing the spectral gap.
More precisely, we have the following rephrasing  of  the Aldous--Fill conjecture.

\begin{conjecture}\rm  The spectral gap (algebraic connectivity) of a    $d$-regular graph on $n$ vertices is at least $(1+o(1))\frac{2d\pi^2}{3n^2}$, and the bound is attained at least for one value of $d$.
\end{conjecture}

It is worth mentioning that in \cite{actt}, it is proved that the maximum relaxation time for the random walk on a   graph on $n$ vertices is $(1 +o(1))\frac{n^3}{54}$, settling another conjecture by Aldous and Fill \cite[p.~216]{aldous2002reversible}.

As usual, we denote the  minimum degree of a graph $H$ by $\delta=\delta(H)$ and its diameter by $\dm(H)$.
Let $G$ be a $d$-regular graph and $H$ a graph with  $\delta=d$, both of order $n$.
We say that $G$ is a {\em $\mu$-minimal} $d$-regular graph if $G$ has the smallest $\mu$ among all $d$-regular graphs of order $n$.
Also $H$ is said to be a {\em $\mu$-minimal} graph with $\delta=d$ if $H$ has the smallest $\mu$ among all graphs with $\delta=d$ and order $n$.

Recall that a {\em block} of a graph is a maximal connected subgraph with no cut vertex. The
blocks of a graph fit together in a tree-like structure, called the {\em block-tree} of $G$.
 When $G$ has at least two blocks and its block-tree is a path, we say that $G$ is {\em path-like}. In such a case, $G$ has two pendant blocks, which are called {\em end blocks} of $G$.

\subsection{Structure of $\mu$-minimal graphs}

L.  Babai (see \cite{Guiduli}) made a conjecture that described the structure of  $\mu$-minimal cubic (i.e. $3$-regular) graphs. Guiduli \cite{Guiduli} (see also \cite{GuiduliThesis}) proved that $\mu$-minimal cubic graphs are path-like, built from specific blocks.
The result of Guiduli was improved  later  by Brand, Guiduli, and Imrich \cite{Imrich}. They completely characterized $\mu$-minimal cubic  graphs and confirmed the Babai conjecture. For every even $n$, such a graph is proved to be unique. (Cubic graphs always have even orders.)
 Abdi, Ghorbani and Imrich~\cite{AbGhIm} showed that the algebraic connectivity  of these graphs is $(1+o(1))\frac{2\pi^2}{n^2}$, confirming  the Aldous--Fill conjecture for $d=3$.
 Guiduli \cite[Problem~5.2]{GuiduliThesis} asked for a generalization of the aforementioned result of Brand, Guiduli, and Imrich, namely the characterization of $\mu$-minimal $d$-regular graphs. 
In this direction, Abdi and Ghorbani \cite{AbGh} gave a `near' complete characterization\footnote{In \cite{AbGhIm}, it was conjectured that a $\mu$-minimal quartic graph has the following structure: any middle block is $M_4$ (refer to Figure~\ref{fig:Type}), and each end block is one of the four specified blocks. This conjecture has been nearly proven in \cite{AbGh} by allowing one additional end block.} 
 of $\mu$-minimal quartic  (i.e. $4$-regular) graphs.  Based on that,  they established   the Aldous--Fill conjecture for $d=4$.

 Guiduli and Mohar
 proposed another generalization of the Babai conjecture by considering graphs with $\delta=d$ rather than $d$-regular graphs.
 They put forward  the following two conjectures on the structure of  $\mu$-minimal graphs with $\delta=d$.
\begin{conjecture}[Guiduli and Mohar, see {\cite[p.~87]{GuiduliThesis}}]\label{conjecture:mohar1}
\rm Let  $n\equiv0 \pmod{d+1}$. Then the $\mu$-minimal graph on $n$ vertices with  $\delta=d$  is the graph of  Figure~\ref{fig:Mohar}.
\end{conjecture}
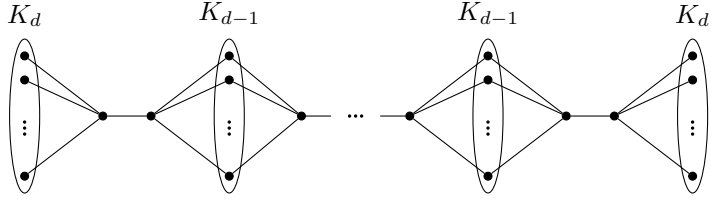
\begin{figure}[t]
\centering
\begin{tikzpicture}[scale=.8]
  \draw (1,0) ellipse (.25 and 1.28);
  \vertex[fill] (1) at (1,1) [] {};
  \vertex[fill] (2) at (1,.6) [] {};
  \vertex[fill] (3) at (1,-1) [] {};
  \vertex[fill] (4) at (2.3,0) [] {};
  \vertex[fill] (5) at (3.1,0) [] {};
  \draw (4.4,0) ellipse (.25 and 1.28);
  \vertex[fill] (6) at (4.4,1) [] {};
  \vertex[fill] (7) at (4.4,.6) [] {};
  \vertex[fill] (8) at (4.4,-1) [] {};
  \vertex[fill] (9) at (5.6,0) [] {};
  \vertex[fill] (12) at (7.4,0) [] {};
  \draw (8.7,0) ellipse (.25 and 1.28);
  \vertex[fill] (13) at (8.7,1) [] {};
  \vertex[fill] (14) at (8.7,.6) [] {};
  \vertex[fill] (15) at (8.7,-1) [] {};
  \vertex[fill] (16) at (10,0) [] {};
  \vertex[fill] (17) at (10.8,0) [] {};
  \draw (12.1,0) ellipse (.25 and 1.28);
  \vertex[fill] (18) at (12.1,1) [] {};
  \vertex[fill] (19) at (12.1,.6) [label=above:\footnotesize{}] {};
  \vertex[fill] (20) at (12.1,-1) [label=below:\footnotesize{}] {};
  \tikzstyle{vertex}=[circle, draw, inner sep=.3pt, minimum size=.3pt]
  \vertex[fill] () at (1,-.1) [label=below:\footnotesize{}] {};
  \vertex[fill] () at (1,-.2) [label=below:\footnotesize{}] {};
  \vertex[fill] () at (1,-.3) [label=below:\footnotesize{}] {};
  \vertex[fill] () at (4.4,-.1) [label=below:\footnotesize{}] {};
  \vertex[fill] () at (4.4,-.2) [label=below:\footnotesize{}] {};
  \vertex[fill] () at (4.4,-.3) [label=below:\footnotesize{}] {};
  \vertex[fill] () at (6.4,0) [label=left:\footnotesize{}] {};
  \vertex[fill] () at (6.5,0) [label=left:\footnotesize{}] {};
  \vertex[fill] () at (6.6,0) [label=left:\footnotesize{}] {};
  \vertex[fill] () at (8.7,-.1) [label=below:\footnotesize{}] {};
  \vertex[fill] () at (8.7,-.2) [label=below:\footnotesize{}] {};
  \vertex[fill] () at (8.7,-.3) [label=below:\footnotesize{}] {};
  \vertex[fill] () at (12.1,-.1) [ ] {};
  \vertex[fill] () at (12.1,-.2) [ ] {};
  \vertex[fill] () at (12.1,-.3) [label=below:\footnotesize{}] {};
  \tikzstyle{vertex}=[circle, draw, inner sep=0pt, minimum size=0pt]
  \vertex[] () at (1,1.3) [label=above:\footnotesize{$K_{d}$}] {};
  \vertex[] () at (4.4,1.3) [label=above:\footnotesize{$K_{d-1}$}] {};
  \vertex[] (10) at (6.1,0) [label=left:\footnotesize{}] {};
  \vertex[] (11) at (6.9,0) [label=left:\footnotesize{}] {};
  \vertex[] () at (8.7,1.3) [label=above:\footnotesize{$K_{d-1}$}] {};
  \vertex[] () at (12.1,1.3) [label=above:\footnotesize{$K_{d}$}] {};
  \path
	        (1) edge (4)
	        (2) edge (4)
	        (3) edge (4)
	        (4) edge (5)
	        (6) edge (5)
	        (7) edge (5)
	        (8) edge (5)
	        (6) edge (9)
	        (7) edge (9)
	        (8) edge (9)
	        (9) edge (10)
	        (11) edge (12)
	        (12) edge (13)
	        (12) edge (14)
	        (12) edge (15)
	        (16) edge (13)
	        (16) edge (14)
	        (16) edge (15)
	        (16) edge (17)
	        (17) edge (18)
	        (17) edge (19)
	        (17) edge (20)  ;
\end{tikzpicture}
\caption{The  conjectured $\mu$-minimal graph with  $\delta=d$ and order  $n\equiv0 \pmod{d+1}$.  Here $K_l$ is the complete graph of order $l$.}\label{fig:Mohar}
\end{figure}
For general $n$, they  conjectured that $\mu$-minimal graphs have almost the same structure:
\begin{conjecture}[Guiduli and Mohar, see {\cite[p.~88]{GuiduliThesis}}]\label{conjecture:mohar2}
\rm  Let $G$ be a $\mu$-minimal graph with $\delta=d$.  Then $G$ is path-like, and except for some blocks near each end, the graph has the same structure as Figure~\ref{fig:Mohar}.
\end{conjecture}
 A more precise phrasing of Conjecture~\ref{conjecture:mohar2} is that for every integer $d$, there exist constants $C_1$ and $C_2$ such that any $\mu$-minimal graph with $\delta=d$ and order at least $C_1$ is path-like and except for a limited number of blocks positioned at either end of the path representing the block-tree of $G$ and containing at most $C_2$ vertices in total, the remaining blocks exhibit the structure of Figure~\ref{fig:Mohar}.

\begin{figure}
	\centering
	{\begin{tikzpicture}[scale=.8]
			\draw (-1.9+.2,0) ellipse (.25 and 1.28);
			\vertex[fill] (05) at (-.6+.2,.4) [ ] {};
			\vertex[fill] (04) at (-.6+.2,-.4) [ ] {};
			\vertex[fill] (03) at (-1.9+.2,1) [ ] {};  
			\vertex[fill] (02) at (-1.9+.2,.6) [ ] {};
			\vertex[fill] (01) at (-1.9+.2,-1) [ ] {};
			\vertex[fill] (0) at (-3+.2,0) [ ] {};
			\draw (1+.2,0) ellipse (.25 and 1.28);
			\vertex[fill] (1) at (-.1+.2,0) [ ] {};
			\vertex[fill] (5) at (2.3+.2,.4) [ ] {};
			\vertex[fill] (55) at (2.3+.2,-.4) [ ] {};
			\vertex[fill] (555) at (2.8+.2,0) [ ] {};
			\vertex[fill] (2) at (1+.2,1) [ ] {};
			\vertex[fill] (3) at (1+.2,.6) [ ] {};
			\vertex[fill] (4) at (1+.2,-1) [ ] {};
	
			\draw (5.7-.2,0) ellipse (.25 and 1.28); 
			\vertex[fill] (6) at (4.6-.25,0) [ ] {};
			\vertex[fill] (07) at (5.7-.25,1) [ ] {};
			\vertex[fill] (8) at (5.7-.25,.6) [ ] {};
			\vertex[fill] (9) at (5.7-.325,-1) [ ] {};
			\vertex[fill] (10) at (7-.25,.4) [ ] {};
			\vertex[fill] (11) at (7-.25,-.4) [ ] {};
			\vertex[fill] (12) at (7.5-.25,0) [ ] {};
			\draw (8.6-.2,0) ellipse (.25 and 1.28);
			\vertex[fill] (13) at (8.6-.25,1) [ ] {};
			\vertex[fill] (14) at (8.6-.25,.6) [ ] {}; 
			\vertex[fill] (15) at (8.6-.25,-1) [ ] {};
			\vertex[fill] (16) at (9.9-.25,.4) [ ] {}; 
			\vertex[fill] (17) at (9.9-.25,-.4) [ ] {};
			\vertex[fill] (18) at (10.4-.25,0) [ ] {};
			\tikzstyle{vertex}=[circle, draw, inner sep=0pt, minimum size=0pt]  
			\vertex[fill] (a) at (10.7-.25,.27) [ ] {}; 
			\vertex[fill] (b) at (10.7-.25,.15) [ ] {}; 
			\vertex[fill] (c) at (10.7-.25,.-.27) [ ] {}; 
			\vertex[fill] (055) at (-3.3+.2,.27) [ ] {};
			\vertex[fill] (0555) at (-3.3+.2,-.27) [ ] {};
			\vertex[fill] (7) at (3.1+.2,.27) [ ] {};  
			\vertex[fill] (7a) at (3.1+.2,.15) [ ] {}; 
			\vertex[fill] (77) at (3.1+.2,-.27) [ ] {};
			\vertex[fill] (777) at (4.3-.25,.27) [ ] {};  
			\vertex[fill] (7777) at (4.3-.25,-.27) [ ] {};
			\vertex[fill] () at (1+.2,1.3) [label=above:\footnotesize{$K_{d-2}$}] {};
			\vertex[fill] () at (-1.9+.2,1.3) [label=above:\footnotesize{$K_{d-2}$}] {};
			\vertex[fill] () at (5.7+.2,1.3) [label=above:\footnotesize{$K_{d-2}$}] {};
			\vertex[fill] () at (8.6+.2,1.3) [label=above:\footnotesize{$K_{d-2}$}] {};
			\tikzstyle{vertex}=[circle, draw, inner sep=.1pt, minimum size=.1pt]
			\vertex[fill] () at (3.1+.2,-.01) [ ] {};  
			\vertex[fill] () at (3.1+.2,-.06) [ ] {}; 
			\vertex[fill] () at (3.1+.2,-.11) [ ] {};
			\vertex[fill] () at (10.7-.25,-.01) [ ] {};  
			\vertex[fill] () at (10.7-.25,-.06) [ ] {}; 
			\vertex[fill] () at (10.7-.25,-.11) [ ] {};
			\tikzstyle{vertex}=[circle, draw, inner sep=.3pt, minimum size=.3pt]
			\vertex[fill] () at (1+.2,-.060) [ ] {};
			\vertex[fill] () at (1+.2,-.160) [ ] {};
			\vertex[fill] () at (1+.2,-.260) [ ] {};
			\vertex[fill] () at (-1.9+.2,-.060) [ ] {};
			\vertex[fill] () at (-1.9+.2,-.160) [ ] {};
			\vertex[fill] () at (-1.9+.2,-.260) [ ] {};
			\vertex[fill] () at (5.7-.25,-.060) [ ] {};
			\vertex[fill] () at (5.7-.25,-.160) [ ] {};
			\vertex[fill] () at (5.7-.25,-.260) [ ] {};
			\vertex[fill] () at (8.6-.25,-.060) [ ] {};
			\vertex[fill] () at (8.6-.25,-.160) [ ] {};
			\vertex[fill] () at (8.6-.25,-.260) [ ] {};
		    	\vertex[fill] () at (3.6,0) [ ] {};  
			    \vertex[fill] () at (3.7,0) [ ] {};  
			    \vertex[fill] () at (3.8,0) [ ] {};  
			\path
			(18) edge (a)
			(18) edge (b)
			(18) edge (c)
			(0) edge (01)
			(0) edge (055)
			(0) edge (0555)
			(0) edge (02)
			(0) edge (03)
			(04) edge (01)
			(04) edge (02)
			(04) edge (03)
			(05) edge (01)
			(05) edge (02)
			(05) edge (03)
			(04) edge (1)
			(05) edge (1) 
			(05) edge (04)  
			(1) edge (2) 
			(1) edge (3) 
			(1) edge (4) 
			(5) edge (2)
			(5) edge (3)
			(5) edge (4)
			(55) edge (2)
			(55) edge (3)
			(55) edge (4)
			(55) edge (5)
			(555) edge (5)
			(555) edge (55)
			(555) edge (7)
			(555) edge (7a)
			(555) edge (77)
			(6) edge (777)
			(6) edge (7777)
			(6) edge (07)
			(6) edge (8)
			(6) edge (9)
			(10) edge (07)
			(10) edge (8)
			(10) edge (9)
			(11) edge (07)
			(11) edge (8)
			(11) edge (9)
			(11) edge (10)
			(12) edge (10)
			(12) edge (11)
			(12) edge (14)
			(12) edge (15)
			(12) edge (13)
			(16) edge (14)
			(16) edge (15)
			(16) edge (13)
			(17) edge (14)
			(17) edge (15)
			(17) edge (13)
			(17) edge (16)
			(17) edge (18)
			(16) edge (18) ;                     
	\end{tikzpicture}}\caption{The conjectured structure (of middle blocks) of $\mu$-minimal $d$-regular graphs for even $d$.}\label{fig:d-even}
\end{figure}
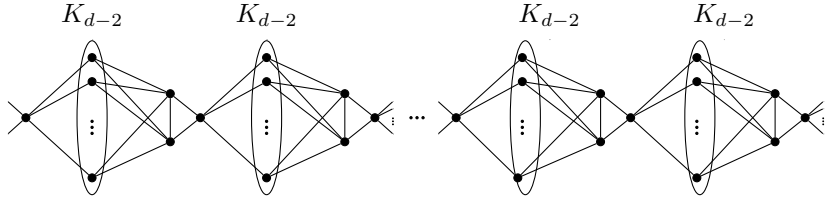

Returning to regular graphs, when $d$ is odd, it is possible to construct $d$-regular graphs with the structure outlined in Conjecture~\ref{conjecture:mohar2} by selecting suitable end blocks. These graphs emerge as natural candidates for $\mu$-minimal $d$-regular graphs. However, for even~$d$, such a construction is not applicable, primarily because regular graphs with even degrees have no bridges.
In this case, we conjecture that $\mu$-minimal regular graphs should exhibit a different structure, as illustrated in Figure~\ref{fig:d-even}. To summarize, we have the following conjecture:

\begin{conjecture}\label{conjecture:ours}
\rm For every integer $d\ge3$, there exist constants $C_1$ and $C_2$ such that any $\mu$-minimal $d$-regular graph $G$ of order at least $C_1$ is path-like and except for a limited number of blocks at either end containing at most $C_2$ vertices in total, $G$ has the same structure as Figure~\ref{fig:Mohar} for odd $d$ and  as Figure~\ref{fig:d-even} for even $d$.
\end{conjecture}

As one of the main results of this paper, we prove that:
\begin{restatable}[]{theorem}{firsttheorem}\label{thm:ours=aldous}
Conjecture~\ref{conjecture:ours} implies  the Aldous--Fill conjecture.	
\end{restatable}
This in particular means that the Guiduli--Mohar conjecture (Conjecture~\ref{conjecture:mohar2}) implies  the Aldous--Fill conjecture for odd $d$.

\subsection{Graphs with maximum diameter}
The maximum diameter of $d$-regular graphs (or those with $\delta=d$) of order $n$ is about $3n/(d+1)$ (see Theorems~\ref{thm:diameter}
and \ref{thm:diameterReg} below). 
The conjectured $\mu$-minimal graphs of Conjectures~\ref{conjecture:mohar2} and \ref{conjecture:ours}  achieve  this maximum diameter.
This phenomenon has been already noted in the literature.
According to Godsil and Royle \cite[p.~289]{Godsil}:
``It has been noted empirically that $\mu(G)$ seems to give a fairly natural measure of 	the `shape' of a graph. Graphs with small values of $\mu(G)$ tend to be elongated graphs of large diameter with bridges.''
Guiduli \cite[p. 46]{GuiduliThesis} showed that  the unique $\mu$-minimal cubic graph 
has the maximum diameter among  cubic graphs of order $n$.
For $n\equiv2\pmod4$, the graph is also the unique one with maximum diameter.
This is not the case for  $n\equiv0\pmod{4}$, 
where there are
$\lfloor (n-4)/8\rfloor$ graphs with the maximum diameter.
Hence he  posed the following problem:
 \begin{problem}[{Guiduli \cite[p. 87]{GuiduliThesis}}]\label{problem1}
\rm Is it true that the cubic graphs  with maximal diameter have algebraic connectivity smaller than all others?
\end{problem}

We show that the answer to this problem in its general form, i.e., for $d$-regular graphs for every $d\ge3$, is negative.
We then consider the asymptotic variant of Problem~\ref{problem1}. In this regard, we establish that $d$-regular graphs for $d\ge5$, as well as graphs with $\delta=d$ for $d\ge4$ with asymptotically maximum diameter (that is $(1+o(1))\frac{3n}{d+1}$) do not necessarily exhibit the asymptotically smallest $\mu$.
For $3$- and $4$-regular graphs, however, we show that a weaker version of the asymptotic problem holds. We conjecture that the converse of the asymptotic variant of Problem~\ref{problem1} is true.
The above results rely on our understanding of the structure as well as optimal estimation of the algebraic connectivity of graphs with diameter $\frac{3n}{d+1}+O(1)$.
Based on that, we also conclude the following theorem which, in particular, implies the Aldous--Fill conjecture for graphs with diameter $\frac{3n}{d+1}+O(1)$.

\begin{restatable}[]{theorem}{secondtheorem} \label{thm:muDiamter3nd+1}
Given $d\ge3$, among graphs with diameter
	$\frac{3n}{d+1}+O(1)$,
	the minimum algebraic connectivity
\begin{itemize}
	\item[\rm(i)] for  graphs with $\delta=d$ is $(1+o(1))\frac{(d-1)\pi^2}{n^2}$,
	\item[\rm(ii)] for  $d$-regular graphs is $(1+o(1))\frac{(d-1)\pi^2}{n^2}$ if $d$ is odd  and $(1+o(1))\frac{2(d-2)\pi^2}{n^2}$ if $d$ is even.
\end{itemize}
In particular, the maximum relaxation time among all
regular graphs  with diameter 	$\frac{3n}{d+1}+O(1)$ is $(1+o(1)) \frac{3n^2}{2\pi^2}$ and is achieved by cubic graphs.
\end{restatable}

 The rest of the paper is organized as follows.
In Section~\ref{sec:o(1/n)}, we establish some properties of graphs with $\mu=o(1/n)$. These results are crucial for our asymptotic arguments.
Section~\ref{sec:MaxDiam} is devoted to nearly-maximum diameter graphs. We give a characterization of such graphs and estimate their algebraic connectivity.  The proof of  Theorems~\ref{thm:ours=aldous} and \ref{thm:muDiamter3nd+1}  will be given in Section~\ref{sec:proofs}. In Section~\ref{sec:MaxMin}, we answer Problem~\ref{problem1} and its generalization to $d$-regular as well as graphs with $\delta=d$ and go through their asymptotic formulations.

\section{Graphs with algebraic connectivity  $o(1/n)$}\label{sec:o(1/n)}
 An eigenvector corresponding to $\mu(G)$ is known as a {\em Fiedler vector}.
In this section, we  extract some facts on the magnitude of the components of a unit  Fiedler vector of a graph of order $n$ and  $\mu = o(1/n)$.
 Then we establish that in such a graph, a perturbation of size $O(1)$ does not change the order of $\mu$.
 These  results will be used in the next sections.

Recall that for a graph $G$ of order $n$ with Laplacian matrix $L(G)$ and $\x\in \mathbb{R}^n$, the quantity $\frac {\x^\top L(G)\x}{ \|\x\|^2}$ is called a {\em Rayleigh quotient}.
   It is well known that
\begin{equation}\label{eq:ReqR}
\mu(G)=\min_{\x\ne\bf0,\,\x\perp\bf1}\frac{\x^\top L(G)\x}{\|\x\|^2},
\end{equation}
where $\bf1$ is the all-$1$ vector.

The quantity $\x^\top L(G)\x$ with $\x=(x_1,\ldots,x_n)^\top$  can be expressed in the following useful manner:
\begin{equation}\label{eq:yL(G)y}
\x^\top L(G)\x = \sum_{ij\in E(G)}(x_i-x_j)^2,
\end{equation}
where $E(G)$ is the edge set of $G$.
Note that if $\x$ is  an eigenvector  for $L(G)$ corresponding to $\mu$,  then  for any vertex $i$ with degree $d_i$,
\begin{equation}\label{eq:eigenequation}
\mu x_i=d_ix_i-\sum_{j:\,ij\in E(G)}x_j.
\end{equation}
We refer to \eqref{eq:eigenequation} as the {\em eigen-equation}.
This also can be written as
\begin{equation*}
\mu x_i=\sum_{j:\,ij\in E(G)}(x_i-x_j).
\end{equation*}

The following lemma allows us to extend \eqref{eq:ReqR} to vectors that are not necessarily orthogonal to $\bf1$. The notation $\langle\cdot,\cdot\rangle$ as usual denotes the standard inner product of real vectors.
 \begin{lemma}\label{lem:remarkdelta}
		Let $G$ be a graph of order $n$ and $\x$ be a vector of length $n$ which is not a multiple of $\1$ and
	  $\|\x\|$ is greater than a positive constant.
	 If $\langle\x,\1\rangle=o(\sqrt n)$, then
		$$\mu(G)\le(1+o(1))\frac{\x^\top L(G)\x}{\|\x\|^2}.$$
\end{lemma}
\begin{proof}
 Let $\e:=\langle\x,\1\rangle$ and
	$\y=\x-\frac\e n\1$. Then $\y\perp\1$, and
	$$\|\y\|^2=\|\x\|^2-\frac{2\e} n \langle\x,\1\rangle+\frac{\e^2}n=\|\x\|^2-\frac{\e^2}n.$$
	 Furthermore, since $L(G)\1=\bf0$, we have $\y^\top L(G)\y=\x^\top L(G)\x$.
	Since $\y\perp\1$, $\mu(G)\le\frac{\y^\top L(G)\y}{\|\y\|^2}$, and thus
	$$\mu(G)\le\frac{\x^\top L(G)\x}{\|\x\|^2-\frac{\e^2}n}.$$
	The right-hand side is $(1+o(1))\frac{\x^\top L(G)\x}{\|\x\|^2}$
as $\|\x\|$ is bounded away from zero and $\e^2/n=o(1)$.		
	\end{proof}

The next lemma illustrates that if $\mu=o(1/n)$, then the components of a unit Fiedler vector tend to $0$ as $n$ grows.

\begin{lemma}\label{lem:o(1)}
  Let $G$ be a graph with $n$ vertices and algebraic connectivity $\mu=o(1/n)$.
 If $\x$ is a unit eigenvector corresponding to $\mu$, then each component of $\x$  is $o(1)$.
\end{lemma}
\begin{proof}With no loss of generality assume that $x_1$ and $x_\ell$ (corresponding to the vertices $v_1$ and $v_\ell$) are the components of $\x$ with the maximum and minimum absolute values, respectively. It suffices to show that $x_1=o(1)$. As $\|\x\|=1$, it is clear that $x_\ell=o(1)$.
There is a path in $G$ between $v_1$ and $v_\ell$.
 With no loss of generality we may assume that $v_1v_2\ldots v_\ell$ is that path.
We have
\begin{align*}
(x_1-x_\ell)^2&=\left(\sum_{i=1}^{\ell-1}(x_i-x_{i+1})\right)^2\\
&\le	(\ell-1)\sum_{i=1}^{\ell-1}(x_i-x_{i+1})^2\\
&\le	(\ell-1)\sum_{ij\in E(G)}(x_i-x_j)^2\\
&=(\ell-1)\mu
\\& \leq n\, o(1/n)\\
&=o(1).	
\end{align*}
This implies that $x_1=o(1)$.
\end{proof}	

\begin{lemma}\label{lem:muedge}
  Let $G$ be a  graph of order $n$ and algebraic connectivity $\mu=o(1/n)$.
Let $\x$ be a unit Fiedler vector of $G$.
    If  $x_r$ and $x_s$ are two components of $\x$ corresponding to vertices at distance $O(1)$,  then $(x_r-x_s)^2=o(\mu)$.
  \end{lemma}
  \begin{proof}
 Let $x_r$ and $x_s$ represent two vertices of distance $t$.   	
 First, assume that $t=1$.
  	With no loss of generality we can assume that $x_r>x_s$.
  	Let $R$ be the set of vertices whose components in $\x$ are greater than or equal to $x_r$.
  	Then it is clear that for $i\in R$ and $j\in S:= V(G)\setminus R$ one has $x_i>x_j$.
  	By applying the eigen-equation %\eqref{eq:eigenequation} 
  	to the vertices of $R$, we have
  	$$\mu \sum_{i\in R}x_i=\sum_{i\in R,\, j\in V(G)\atop i\sim j}(x_i-x_j)=\sum_{i\in R,\, j\in S\atop i\sim j}(x_i-x_j).$$
 (The edges with both endpoints in $R$ contribute $0$ to the middle sum.)
  In the right-hand sum,	every term  is positive and additionally one of its term is $(x_r-x_s)$. It follows that
  $$(x_r-x_s)^2 \le \mu^2\left(\sum_{i\in R}x_i\right)^2\le \mu^2|R|\sum_{i\in R}x_i^2\le\mu^2n=o(\mu).$$
 
 Now, suppose that $t>1$. So we can assume that $x_r=x'_0,x'_1,\ldots,x'_t=x_s$ are the components of $\x$ corresponding to the vertices of a path of length $t$. Then
  $$(x_r-x_s)^2=\left(\sum_{i=1}^t(x'_i-x'_{i-1})\right)^2\le t\sum_{i=1}^t(x'_i-x'_{i-1})^2\le  t^2 o(\mu).$$
 The result now follows since $t=O(1)$.
  \end{proof}

In the final result of this section,  we demonstrate that for a graph with a small enough $\mu$, a perturbation of size $O(1)$ changes its algebraic connectivity only by $o(\mu)$.

\begin{theorem} \label{thm:H}
	Let $G$ be a graph of order $n$ and $\mu(G)=o(1/n)$. Let $H$ be another graph and
	$G'$ be  a connected graph obtained from $G$ by connecting some vertices of $H$ to the vertices in $S\subseteq V(G)$.
	If $S$ and  $H$ are both of order $O(1)$  and 
	 the distance of any pair of vertices of $S$ in $G$ is also $O(1)$, then $\mu(G')=(1+o(1))\mu(G)$.
\end{theorem}
\begin{proof}
Let  $\x$ be a unit Fiedler vector of $G$ and $\mu=\mu(G)$. Let $x_0$ be a component of $\x$ corresponding to some fixed vertex of $S$.
As  the distance of any pair of vertices of $S$ in $G$  is $O(1)$, by Lemma~\ref{lem:muedge}, 
\begin{equation}\label{eq:(x_s-x)^2}
(x_0-x_s)^2=o(\mu)\quad\hbox{for any component $x_s$ of $\x$ corresponding to a vertex in $S$.}
\end{equation}
 Let $H$ have $k$ vertices. We extend $\x$ to a vector $\x'$ of length $n+k$ on $G'$ as follows: on $H$, all the components of $\x'$ are equal to $x_0$, and on   the remaining vertices, $\x'$ agrees with $\x$.
  So,  by considering \eqref{eq:yL(G)y} and
  \eqref{eq:(x_s-x)^2}, $\x'^\top L(G')\x'=(1+o(1))\mu$. We have $\langle\x',\1_{n+k}\rangle=\langle\x,\1_n\rangle+k x_0=k x_0$ which is $o(1)$
  by Lemma~\ref{lem:o(1)}. Similarly $\|\x'\|^2=1+o(1)$.  Thus by Lemma~\ref{lem:remarkdelta},
$$\mu(G')\le(1+o(1))\frac{\x'^\top L(G')\x'}{\|\x'\|^2}=(1+o(1))\mu.$$

To establish the reverse inequality, let $\y'$ be a unit Fiedler vector of $G'$ and $\y$ be the restriction of $\y'$ to $G$.
The graph $G'$ has $n'=n+k$ vertices.
Since $\mu(G')\leq (1+o(1)) \mu$, we have $\mu(G')=o(1/n')$.
In view of \eqref{eq:yL(G)y} and by Lemma~\ref{lem:muedge},
all the terms  $(x_i-x_j)^2$ appearing in $\y'^\top L(G')\y'-\y^\top L(G)\y$ are $o(\mu)$ and thus
 $\y^\top L(G)\y=(1+o(1))\mu(G')$.
 On the other hand,  we see that $\langle\y,\1_{n}\rangle=\langle\y',\1_{n'}\rangle-\sum_{v_i\in V(H)}y'_i=-\sum_{v_i\in V(H)}y'_i$ which is
 $o(1)$  by Lemma \ref{lem:o(1)}. Also $\|\y\|^{2}=\|\y'\|^{2}- \sum_{v_i \in V(H)}{y'_i}^2=1+o(1)$.
Therefore, by Lemma~\ref{lem:remarkdelta},
$$\mu\le(1+o(1))\frac{\y^\top L(G)\y}{\|\y\|^2}=(1+o(1))\mu(G'),$$
which completes the proof.
\end{proof}

\section{Nearly maximum-diameter graphs }\label{sec:MaxDiam}
We know that (\cite{Caccetta}, see also Theorem~\ref{thm:diameter} below) for $d\geq3$ and  $n\geq2d+4$, the maximum diameter of a  graph with order $n$ and $\delta=d$ is
$3\lfloor\frac{n}{d+1}\rfloor-\ell$ for some $\ell\in\{1,2,3\}$.
In this section, we investigate graphs  with order $n$,  $\delta=d$, and  diameter
$3n/(d+1)+O(1)$. We  determine their structure and estimate their algebraic connectivity. From these results,  we deduce Theorem~\ref{thm:muDiamter3nd+1} in the next section.

\subsection{The structure}
  Before proceeding, a definition and some notation are in order.
 A partition $\Pi= \lbrace C_1,  \ldots , C_{m} \rbrace$ of $V(G)$
is called an {\em equitable partition} for $G$ if for every pair of (not necessarily distinct) indices $i, j \in \lbrace 1,  \ldots , m\rbrace$, there is a non-negative integer $q_{ij}$ such that each vertex $v$ in the cell $C_i$ has exactly $q_{ij}$ neighbors in the cell $C_j$, regardless of the choice of $v$. 
The {\em sequential join} of  vertex-disjoint graphs $G_1,G_2, \ldots, G_k$, denoted by  $G_1+ G_2+ \cdots+G_k$, is obtained from the union $G_1\cup G_2 \cup \cdots \cup G_k$
by adding edges joining each vertex of $G_i$ with each vertex of $G_{i+1}$ for $i=1,\ldots,k-1$.
We use the notation $\G(a,b,c;m)$ to denote the sequential join of the sequence of $3m$ complete graphs
$K_a,K_b,K_c, K_a,K_b,K_c,\ldots, K_a,K_b,K_c$.
%The graph $(K_a+K_b+K_c)_m$ is also denoted by } 
So this graph has $(a+b+c)m$ vertices.
As an instance, the graph $\G(2,3,4;3)$ is illustrated in Figure~\ref{fig:G(2,3,4,3)}.
Any of the cliques  $K_a$, $K_b$ or $K_c$ (whose  vertices are  drawn vertically above each other in Figure~\ref{fig:G(2,3,4,3)}) in $\G(a,b,c;m)$ will be referred to as a cell.  Such  cliques are   in fact  the cells of the `natural' equitable partition of  the graph.

\begin{figure}[h!]
%	\captionsetup[subfigure]{labelformat=empty}
	\centering
	\begin{tikzpicture}[scale=2]
		\vertex[fill] (1) at (0,.1) [] {};
		\vertex[fill] (2) at (0,-.1) [] {};
		\vertex[fill] (3) at (.5,.2) [] {};
		\vertex[fill] (4) at (.5,0) [] {};
		\vertex[fill] (5) at (.5,-.2) [] {};
		\vertex[fill] (7) at (1,.1) [] {};
		\vertex[fill] (6) at (1,.3) [] {};
		\vertex[fill] (8) at (1,-.1) [] {};
		\vertex[fill] (9) at (1,-.3) [] {};
		\vertex[fill] (a1) at (1.5,.1) [] {};
		\vertex[fill] (a2) at (1.5,-.1) [] {};
		\vertex[fill] (a3) at (2,.2) [] {};
		\vertex[fill] (a4) at (2,0) [] {};
		\vertex[fill] (a5) at (2,-.2) [] {};
		\vertex[fill] (a7) at (2.5,.1) [] {};
		\vertex[fill] (a6) at (2.5,.3) [] {};
		\vertex[fill] (a8) at (2.5,-.1) [] {};
		\vertex[fill] (a9) at (2.5,-.3) [] {};
		\vertex[fill] (b1) at (3,.1) [] {};
		\vertex[fill] (b2) at (3,-.1) [] {};
		\vertex[fill] (b3) at (3.5,.2) [] {};
		\vertex[fill] (b4) at (3.5,0) [] {};
		\vertex[fill] (b5) at (3.5,-.2) [] {};
		\vertex[fill] (b7) at (4,.1) [] {};
		\vertex[fill] (b6) at (4,.3) [] {};
		\vertex[fill] (b8) at (4,-.1) [] {};
		\vertex[fill] (b9) at (4,-.3) [] {};
		\path
		(1) edge (2)
		(1) edge (3)
		(1) edge (4)
		(1) edge (5)
		(2) edge (3)
		(2) edge (4)
		(2) edge (5)
		(3) edge (4)
		(4) edge (5)
		(7) edge (2)
		(3) edge[bend right=40] (5)
		(3) edge (6)
		(3) edge (7)
		(3) edge (8)
		(3) edge (9)
		(4) edge (6)
		(4) edge (7)
		(4) edge (8)
		(4) edge (9)
		(5) edge (6)
		(5) edge (7)
		(5) edge (8)
		(5) edge (9)
		(6) edge (7)
		(7) edge (8)
		(8) edge (9)
		(6) edge[bend left=40] (8)
		(6) edge[bend left=60] (9)
		(7) edge[bend left=40] (9)
		(6) edge (a1)
		(6) edge (a2)
		(7) edge (a1)
		(7) edge (a2)
		(8) edge (a1)
		(8) edge (a2)
		(9) edge (a1)
		(9) edge (a2)
		(a1) edge (a2)
		(a1) edge (a3)
		(a1) edge (a4)
		(a1) edge (a5)
		(a2) edge (a3)
		(a2) edge (a4)
		(a2) edge (a5)
		(a3) edge (a4)
		(a4) edge (a5)
		(a7) edge (a2)
		(a3) edge[bend right=40] (a5)
		(a3) edge (a6)
		(a3) edge (a7)
		(a3) edge (a8)
		(a3) edge (a9)
		(a4) edge (a6)
		(a4) edge (a7)
		(a4) edge (a8)
		(a4) edge (a9)
		(a5) edge (a6)
		(a5) edge (a7)
		(a5) edge (a8)
		(a5) edge (a9)
		(a6) edge (a7)
		(a7) edge (a8)
		(a8) edge (a9)
		(a6) edge[bend left=40] (a8)
		(a6) edge[bend left=60] (a9)
		(a7) edge[bend left=40] (a9)	
		(a6) edge (b1)
		(a6) edge (b2)
		(a7) edge (b1)
		(a7) edge (b2)
		(a8) edge (b1)
		(a8) edge (b2)
		(a9) edge (b1)
		(a9) edge (b2)	
		(b1) edge (b2)
		(b1) edge (b3)
		(b1) edge (b4)
		(b1) edge (b5)
		(b2) edge (b3)
		(b2) edge (b4)
		(b2) edge (b5)
		(b3) edge (b4)
		(b4) edge (b5)
		(b7) edge (b2)
		(b3) edge[bend right=40] (b5)
		(b3) edge (b6)
		(b3) edge (b7)
		(b3) edge (b8)
		(b3) edge (b9)
		(b4) edge (b6)
		(b4) edge (b7)
		(b4) edge (b8)
		(b4) edge (b9)
		(b5) edge (b6)
		(b5) edge (b7)
		(b5) edge (b8)
		(b5) edge (b9)
		(b6) edge (b7)
		(b7) edge (b8)
		(b8) edge (b9)
		(b6) edge[bend left=40] (b8)
		(b6) edge[bend left=60] (b9)
		(b7) edge[bend left=40] (b9);
	\end{tikzpicture}
	\caption{The graph $\G(2,3,4;3)$.}
	\label{fig:G(2,3,4,3)}
\end{figure}
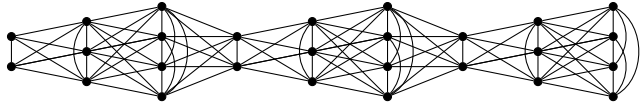

Let $d\ge2$, $t\ge1$,  and  $a_1,b_1,c_1,\ldots,a_t,b_t,c_t$ be positive  integers such that for each $i=1,\ldots,t$ we have $a_i+b_i+c_i=d+1$. Let $\G_i:=\G(a_i,b_i,c_i;m_i)$  which has $n_i=(d+1)m_i$ vertices.
The graph $\g={\Gamma}_d(a_1,b_1,c_1,\ldots,a_t,b_t,c_t;m_1,\ldots,m_t)$ is a graph  obtained  from $\G_1\cup \G_2 \cup \cdots \cup \G_t$ by adding edges joining every vertex of the last cell of $\G_i$ to every vertex of the first cell of  $\G_{i+1}$ for  $i=1, \ldots, t-1$.
 We allow $n=(d+1)(m_1+\cdots+m_t)$
to grow.
In $\g$, every three consecutive cells have $d+1$ vertices, 
except for the triples containing the last cell of $\G_i$ and the first cell of   $\G_{i+1}$.
 So all but at most $a_1+c_1+\cdots+a_t+c_t\le td$ vertices have degree  $d$.
 Since $d$ and $t$ are fixed and $n$ can grow, almost all vertices of $\Gamma$ have degree $d$. Also, $\dm(\g)=3(m_1+\cdots+m_t)-1=3n/(d+1)-1$
which is the maximum diameter of a $d$-regular graph (see Theorem~\ref{thm:diameterReg}).

Finally, we define a family of graphs, namely ${\mathscr F}_{n,d,C}$, which, as we shall prove, characterizes  nearly maximum-diameter graphs with $\delta=d$. 
\begin{definition}\rm
Given positive integers $n,d$ and a constant $C$,
a  graph  $\cal G$ belongs to ${\mathscr F}_{n,d,C}$  if:
\begin{itemize}
\item[(i)] there exist positive integers $t\le C$, $m_1,\ldots,m_t$, and  $a_1,b_1,c_1,\ldots,a_t,b_t,c_t$,
and graphs $H_0,\ldots,H_t$ with $\sum_{i=0}^{t} |V(H_i)|\le C$, 
 such that $a_i+b_i+c_i=d+1$ for  $i=1,\ldots,t$, and  $\sum_{i=1}^t m_i(d+1)+\sum_{i=0}^{t} |V(H_i)|=n$, 
\item[(ii)] $\G$ is connected and obtained form
 $H_0\cup \G_1\cup H_1\cup \G_2 \cup  \cdots \cup H_{t-1}\cup \G_t \cup H_t $, where $\G_i:=\G(a_i,b_i,c_i;m_i)$,
by connecting arbitrary vertices from the first (resp. last) cell of  $\G_{i}$  to arbitrary vertices of $H_{i-1}$ (resp. $H_{i}$). 
\end{itemize}
The graphs $\G_1,\ldots,\G_t$ are called {\em major subgraphs} of $\G$.
\end{definition}

 \begin{theorem}\label{thm:diam} Let $G$ be a graph of order $n$ and $\delta=d$. If $\dm(G)=\frac{3n}{d+1}+O(1)$, then for some constant $C$,
 the graph	$G$ belongs to the family $\mathscr{F}_{n,d,C}$.
\end{theorem}	
\begin{proof} 
 Let $\dm(G)=\ell$, so 
$\ell\ge\frac{3n}{d+1}-c$ for some constant $c\ge1$.
 Consider a  distance-partition $\{P_0, \ldots , P_\ell\}$ of $G$  from a vertex that
is on some longest path, with $p_i=|P_i|$.
Since  $\delta=d$ and the neighbors of a vertex in  $P_{i+1}$ lie in  $P_i\cup P_{i+1}\cup P_{i+2}$, we have $q_i:=p_i+p_{i+1}+p_{i+2}\geq d+1$.
Each vertex of $G$ has a contribution of at most $3$ to the sum $q_0+\cdots+q_{\ell-2}$. It follows that 
\begin{equation}\label{eq:ell+c}
(\ell-1)(d+1)\le q_0+\cdots+q_{\ell-2}\le3n\le(\ell+c)(d+1).
\end{equation}
Let $J:=\{ j\in\{0,\ldots, \ell-2\}: q_j\ge d+2\}$.
From \eqref{eq:ell+c}, we see that $|J|\le(c+1)(d+1)$.
Let $U:=\{j,j+1,j+2 : j\in J\}$.
We can partition $\{0,\ldots, \ell\}$
as  $U_0,V_1, U_1,\ldots, V_t, U_t$ such that
 each $U_i$ and $V_i$ consist of consecutive integers and 
$U_0,U_1\ldots, U_t$ is a partition of $U$.\footnote{Note that $\{0,1,2\}\subseteq U_0$ because
	$p_1\ge d$ (since the  neighbors of the vertex in $P_0$ lie in $P_1$) and so $q_0\ge d+2$. Similarly, $p_{\ell-1}+p_\ell\ge d+1$, so	$\{\ell-2,\ell-1,\ell\}\subseteq U_t$.}
We may further assume that $|V_i|\equiv0\pmod3$, otherwise we remove the last one or two members of $V_i$ and add them to $U_i$. So, we can suppose that $|V_i|=3m_i$ for some positive integer $m_i$.
Let $H_i$ and $G_i$ be the induced subgraphs of $G$ on $\bigcup_{j\in U_i}P_j$, and  $\bigcup_{j\in V_i}P_j$, respectively. 
 Assume that $V_i=\{r+1,\ldots,r+3m_i\}$
which implies that all the consecutive triples 
in the sequence $p_{r+1},\ldots,p_{r+3m_i}$ sum up to $d+1$. This is only possible when the entire sequence is a repetition of the first three terms.
So, $G_i=\G(a_i,b_i,c_i;m_i)$, where
$a_i=p_{r+1}$, $b_i=p_{r+2}$, $c_i=p_{r+3}$.
Let $C:=2(c+1)(d+1)^2$. We have $$\sum_{i=0}^{t} |V(H_i)|\le |J|d+\sum_{j\in J}q_j\le |J|d+ 3n-(\ell-1-|J|)(d+1)\le C.$$
Therefore, we have established that $G\in{\mathscr F}_{n,d,C}$. 
\end{proof}

\subsection{The algebraic connectivity}

We start by estimating the algebraic connectivity of 
$$\g=\g_d(a_1,b_1,c_1,\ldots,a_t,b_t,c_t;m_1,\ldots,m_t).$$
Note that $d,t$ are fixed and $m_1+\cdots+m_t\to\infty$.
 For this purpose, we first analyze the Fiedler vector of $\g$.

\begin{lemma}[\hspace{1sp}Fiedler \cite{Fiedler}] \label{lem:Fiedler1}
 Let $\y$ be a Fiedler vector of a  graph $G$ and vertex set $V=\lbrace v_1,\ldots, v_n\rbrace$. Let $V_1=\lbrace v_i\in V: y_i\geq 0\rbrace$ and $V_2=\lbrace v_i\in V: y_i \leq 0\rbrace$. Then both the subgraphs induced by $V_1$ and $V_2$ are connected.
\end{lemma}

\begin{lemma}[\hspace{1sp}Fiedler \cite{Fiedler}]\label{lem:Fiedler2}
  Let $\y$ be a Fiedler vector of a  graph $G$. If $y_i > 0$,
then there exists a vertex j such that $i\sim j$ and $y_j < y_i$.
\end{lemma}

Now we can infer some useful properties of the Fielder vector of $\g$.
 
\begin{lemma} \label{lem:sign}
 Let $\y$ be a Fiedler vector of ${\Gamma}_d(a_1,b_1,c_1,\ldots,a_t,b_t,c_t;m_1,\ldots,m_t)$.
 Let $m=\sum_{j=1}^tm_j$ and $\Pi=\{C_1, \ldots , C_{3m}\}$ (numbered consecutively from left to right)  be  an equitable partition of the vertex set ${\Gamma}$ in which each cell $C_i$ is a $K_{a_i}$, $K_{b_i}$, or $K_{c_i}$.
 \begin{itemize}
   \item[\rm(i)] The  components of $\y$ on each cell  of the partition $\Pi$ are equal.
   \item[\rm(ii)]  Let $y_1, \ldots, y_{3m}$ be the values of $\y$ on the cells of $\Pi$. Then  the $y_i$'s form a strictly monotone sequence changing sign once.
 \end{itemize}
  \end{lemma}
 \begin{proof}
By using  the eigen-equation, %\eqref{eq:eigenequation}, 
we observe that the components of $\y$ on each cell $C_i$ are equal.
Lemma~\ref{lem:Fiedler1} allows us to assume that for some $r\ge1$ and $s\ge0$, all $y_1,\ldots, y_r$ are positive, all $y_{r+s+1},\ldots,y_{3m}$ are negative, and 
all other $y_i$'s (if any) are zero.
From Lemma~\ref{lem:Fiedler2}, it follows that $y_1>y_2>\cdots>y_r$. Now consider $-\y$ as a
Fiedler vector of $\Gamma$. Again by Lemma~\ref{lem:Fiedler2}, $-y_{3m}>-y_{3m-1}>\cdots >-y_{r+s+1}$.
 Hence $y_{r+s+1}>y_{r+s+2}>\cdots >y_{3m}$.  Therefore, $y_i$'s satisfy (ii).
\end{proof}

The path-like structure of the graphs $\g$ allows one to `approximate' their Fiedler vectors using the Fiedler vectors of paths.
For this reason, we first recall  what the Fiedler vector of a path is.
\begin{remark}\label{rem:muP_n}\rm
 For $P_n$, the path graph on $n$ vertices, we  know that  $\mu(P_n)=2(1-\cos\left(\frac{\pi}{n}\right))$  (see \cite{fiedler1973algebraic}),
 and by \cite[p.~53]{Spielman}, its Fiedler vector is   $(x_1,\ldots, x_n)^\top$  with
$$ x_i=\cos\left(\frac{(2i-1)\pi}{2n}\right),~~i=1,\ldots,{n}.$$
\end{remark}

We start by establishing an optimal upper bound on $\mu(\g)$.

\begin{theorem}\label{thm:upper}
Let ${\Gamma}={\Gamma}_d(a_1,b_1,c_1,\ldots,a_t,b_t,c_t;m_1,\ldots,m_t)$ have order $n$ and $L=\max \{a_jb_jc_j : j=1,\ldots,t\}$. Then $\mu({\Gamma})\leq (1+o(1))\frac{L\pi^2}{n^2}$.
\end{theorem}
\begin{proof}
	Let $m:=\sum_{j=1}^tm_j$. Then  ${\Gamma}$ has $n=(d+1)m$ vertices.
Let $C_1,\ldots, C_{3m}$ be the cells of the equitable partition $\Pi$.
 %and an equitable partition with $3m$ cells in which each cell is the vertex set of either a $K_{a_j}$, $K_{b_j}$, or $K_{c_j}$.
	For $i=1,\ldots,m$, we set
	\begin{equation}\label{eq:cosx}
		x_i:=\sqrt{\frac2m}\cos\left(\frac{(2i-1)\pi}{2m}\right).
	\end{equation}
We assign $x_i$ to the vertices of the cell $C_{3i-2}$. We then extend it to the cells $C_{3i-1}$ and $C_{3i}$ as follows.
Assume that $C_{3i-2}=K_{a_r}$, $C_{3i-1}=K_{b_r}$, and $C_{3i}=K_{c_r}$ for some $a_r,b_r,c_r$. Then
we assign $x'_i$ and $x''_i$ to the vertices of $C_{3i-1}$ and $C_{3i}$, where
	$$x'_i=\frac{(a_r+b_r)x_i+c_rx_{i+1}}{a_r+b_r+c_r}, \quad
x''_i=\frac{b_rx_i+(a_r+c_r)x_{i+1}}{a_r+b_r+c_r}.$$
Further, we set $x_{m+1}$ to be equal to $x_m$, so that $x'_m=x''_m=x_m$.
	These define a vector, say $\y=(y_1,\ldots,y_n)^\top$, on the vertices of ${\Gamma}$.
For $i=1,\ldots,m-1$, let
$G_i$ be the induced subgraph on the four consecutive cells $C_{3i-2},C_{3i-1},C_{3i},C_{3(i+1)-2}$.
For $C_{3(i+1)-2}$ there are two possibilities: it is either $K_{a_r}$ or $K_{a_{r+1}}$.
First assume that the former is the case. Then
	by the definition of $\y$, we have
\begin{align*}
	\sum_{jk\in E(G_i)}(y_j-y_k)^2&=
	a_rb_r(x_i-x'_i)^2+b_rc_r(x'_i-x''_i)^2+c_ra_r(x''_i-x_{i+1})^2\\
	&=\frac{a_rb_rc_r}{a_r+b_r+c_r}(x_i-x_{i+1})^2.
\end{align*}
If $C_{3(i+1)-2}=K_{a_{r+1}}$, then
%\begin{align}\label{eq:upper1}
$$ \sum_{jk\in E(G_i)}(y_j-y_k)^2=\frac{a_rb_rc_r}{a_r+b_r+c_r}(x_i-x_{i+1})^2+c_r(a_{r+1}-a_r)(x''_i-x_{i+1})^2.$$
 %\end{align}
 We see that the second term in the right-hand side is   $O(1/m^3)$. The number of such terms is $t-1=O(1)$.
Moreover, letting $G_m$ to be the induced subgraph on the cells $C_{3m-2},C_{3m-1},C_{3m}$, we have
%\begin{align}\label{eq:upper2}
$$ \sum_{jk\in E(G_m)}(y_j-y_k)^2=0.$$
 %\end{align}
Note that $E(G_1)\cup\cdots\cup E(G_m)$ gives  a partition of  $E({\Gamma})$. It follows that
	\begin{align}\label{eq:uupper}
		\sum_{jk\in E({\Gamma})}(y_j-y_k)^2\leq\frac{L}{d+1} \sum_{i=1}^{m-1}(x_i-x_{i+1})^2 +
   O\left(\frac1{m^3}\right).
	\end{align}
Next we find a lower bound for $\|\y\|^2$.
Let $D_i:=C_{3i-2}\cup C_{3i-1}\cup C_{3i}$ and $q:=\lfloor m/2\rfloor$. We have $x_1>\cdots>x_q>0$. So for $i=1,\ldots,q-1$, both $x'^2_i$ and $x''^2_i$ are greater than $x^2_{i+1}$. It follows that
$$\sum_{j\in D_i}y_j^2=a_rx_i^2+b_rx'^2_i+c_rx''^2_i\geq (d+1)x_{i+1}^2,\quad\hbox{for}~i=1,\ldots,q-1.$$
We have also  $0\ge x_{q+1}>\cdots>x_m$. So for $i=q+1,\ldots,m$, both $x'^2_i$ and $x''^2_i$ are at least $x^2_i$. It follows that
$$\sum_{j\in D_i}y_j^2\geq (d+1)x_i^2,\quad\hbox{for}~i=q+1,\ldots,m.$$
For $\sum_{j\in D_q}y_j^2$ we take into account the trivial lower bound zero.
As $V({\Gamma})=D_1\cup\cdots\cup D_m$, we come up with $\sum_{j=1}^n {y_j}^2   \geq (d+1) \sum_{i=2}^{m} x_i^2$.
Also $x_1^2=O(1/m)$. It follows that
\begin{equation} \label{eq:upperr}
\sum_{j=1}^n {y_j}^2   \geq (d+1) \sum_{i=1}^{m} x_i^2+O\left(\frac1m\right).
\end{equation}
Our next task is to show that $\langle\y,\1\rangle=o(1)$. We have
$$\sum_{j\in D_i}y_j=a_rx_i+b_rx'_i+c_rx''_i=(a_r+b_r)x_i+c_rx_{i+1}.$$
It follows that
$$\sum_{j\in D_1\cup\cdots\cup D_{m_1}}y_j=(a_1+b_1)x_1+(d+1)(x_2+\cdots+x_{m_1})+c_1x_{m_1+1}.$$
By \eqref{eq:cosx}, $x_j=o(1)$, and thus
$$\sum_{j\in D_1\cup\cdots\cup D_{m_1}}y_j=(d+1)\sum_{j=1}^{m_1}x_j+o(1).$$
Similarly, for $i=1,\ldots,t-1$, we have
$$\sum_{j\in D_{m_i+1}\cup\cdots\cup D_{m_{i+1}}}y_j=(d+1)\sum_{j=m_i+1}^{m_{i+1}}x_j+o(1).$$
Summing up all these equalities, we obtain
$$\sum_{j\in V({\Gamma})}y_j=(d+1)\sum_{j=1}^mx_j+o(1).$$
From \eqref{eq:cosx} we see that $\sum_{i=1}^{m} x_i=0$ and thus $\langle\y,\1\rangle=o(1)$.
Hence by Lemma~\ref{lem:remarkdelta} and by  \eqref{eq:uupper} and \eqref{eq:upperr} it is inferred that
	\begin{align*}
		\mu({\Gamma})	& \leq(1+o(1)) \frac{\y^\top L({\Gamma})\y}{\|\y\|^2}\\
		& \leq (1+o(1)) \frac{\sum_{ij\in E({\Gamma})}(y_i-y_j)^2}{\sum_{i=1}^n {y_i}^2}\\
		&\leq (1+o(1)) \frac{L}{(d+1)^2} \frac{\sum_{i=1}^{m-1}(x_i-x_{i+1})^2}{\sum_{i=1}^{m} x_i^2}.
	\end{align*}
	From \eqref{eq:cosx} and Remark~\ref{rem:muP_n} we have
	$$\frac{\sum_{i=1}^{m-1}(x_i-x_{i+1})^2}{\sum_{i=1}^{m} x_i^2}=\mu(P_{m})=(1+o(1)) \frac{\pi^2}{m^2},$$
which implies that
	$$	\mu({\Gamma})\le(1+o(1))\frac{L\pi^2}{n^2}.$$
\end{proof}

Now we establish a lower bound on $\mu(\g)$, which is somewhat dual to the upper bound of Theorem~\ref{thm:upper}.

\begin{theorem}\label{thm:lower}
Let ${\Gamma}={\Gamma}_d(a_1,b_1,c_1,\ldots,a_t,b_t,c_t;m_1,\ldots,m_t)$ have order $n$ and $\ell=\min \{a_jb_jc_j : j=1,\ldots,t\}$. Then  $\mu({\Gamma})\geq (1+o(1))\frac{\ell\pi^2}{n^2}$.
 \end{theorem}
 \begin{proof}
Let $\y=(y_1, y_2, \ldots, y_n)^\top$ be a unit Fiedler vector of ${\Gamma}$.
This is constant on each cell of ${\Gamma}$. Also let $\x$ be a vector of length $m$ consisting of the components of $\y$ on the cells $C_1,C_4,\ldots,C_{3m-2}$.
Let  $G_i$ be the induced subgraph on the four consecutive cells $C_{3i-2},C_{3i-1},C_{3i},C_{3(i+1)-2}$.
Let  $u$ and $v$  be the components of $\y$ on the two middle cells of $G_i$. Suppose that $C_{3i-2}=K_{a_r}$.
If $C_{3(i+1)-2}=K_{a_r}$, then
$$\sum_{jk\in E(G_i)}(y_j-y_k)^2= a_rb_r(x_i-u)^2+b_rc_r(u-v)^2+c_ra_r(v-x_{i+1})^2.$$
The right-hand side, considered as a function of $u$ and $v$,  is minimized at
$$u=\frac{(a_r+b_r)x_{i}+c_r x_{i+1}}{a_r+b_r+c_r},
\quad \text{and} \quad
v=\frac{b_r x_{i}+(a_r+c_r)x_{i+1}}{a_r+b_r+c_r}.$$
This implies that
$$	\sum_{jk\in E(G_i)}(y_j-y_k)^2\ge\frac{a_rb_rc_r}{a_r+b_r+c_r}(x_{i}-x_{i+1})^2\ge \frac\ell{d+1}(x_{i}-x_{i+1})^2.$$
If $C_{3(i+1)-2}=K_{a_{r+1}}$, then
$$ \sum_{jk\in E(G_i)}(y_j-y_k)^2=a_rb_r(x_i-u)^2+b_rc_r(u-v)^2+c_ra_r(v-x_{i+1})^2+c_r(a_{r+1}-a_r)(v-x_{i+1})^2.$$
From Theorem~\ref{thm:upper}, $\mu=O(1/n^2)$ and so by Lemma~\ref{lem:muedge}, we have
$(v-x_{i+1})^2=o(1/n^2)$. It follows that
$$	\sum_{jk\in E(G_i)}(y_j-y_k)^2\ge \frac\ell{d+1}(x_{i}-x_{i+1})^2+o\left(\frac1{n^2}\right).$$
Moreover, letting $G_m$ to be the induced subgraph on the cells $C_{3m-2},C_{3m-1},C_{3m}$, we have
%\begin{equation}\label{eq:upper2}
$$ \sum_{jk\in E(G_m)}(y_j-y_k)^2=o\left(\frac1{n^2}\right).$$
 It is inferred that
\begin{equation}\label{eq:lowereven1}
\mu=\mu({\Gamma})=\sum_{ij\in E({\Gamma})}(y_i-y_j)^2 \geq (1+o(1)) \frac\ell{d+1} \sum_{i=1}^{m-1}(x_i-x_{i+1})^2.
\end{equation}
Note that the right-hand side of \eqref{eq:lowereven1} is
$\Theta(1/n^2)$, a fact that will be clarified shortly. 
This justifies the elimination of $t$ terms $o(1/n^2)$. 
Let $D_i:=C_{3i-2}\cup C_{3i-1}\cup C_{3i}$. By Lemma~\ref{lem:sign}, $y_1\ge\cdots\ge y_n$ and $y_i$'s change sign once.  The same also holds for  $x_1,\ldots,x_m$.
Let $q$ be the index such that $x_q>0\ge x_{q+1}$. Then for $i=1,\ldots,q$,
$$	\sum_{j\in D_i}y_j^2= a_rx_{i}^2+b_ru^2+c_rv^2\leq (d+1)x_i^2.$$
Then for $i=q+1,\ldots,m-1$,
$$	\sum_{j\in D_i}y_j^2= a_rx_{i}^2+b_ru^2+c_rv^2\leq (d+1)x_{i+1}^2.$$
It follows that
$$\sum_{j=1}^n y_j^2   \leq (d+1) \sum_{i=1}^q x_i^2+  (d+1)\sum_{i=q+2}^m x_i^2 + \sum_{j\in D_m} y_j^2
	 \leq (d+1) \sum_{i=1}^m x_i^2+ (d+1)y_n^2.$$
	 By Lemma~\ref{lem:o(1)}, $y_n^2=o(1)$, and thus
\begin{align} \label{eq:lower}
	\sum_{i=1}^n {y_i}^2   \leq (d+1) \sum_{i=1}^m x_i^2+o(1).
\end{align}
For each $i=1,\ldots,m$ we have $\sum_{j\in D_i}y_j\leq (d+1)x_i$. This implies that
$$0=\sum_{i=1}^n y_i  \leq (d+1) \sum_{i=1}^m x_i.$$
On the other hand, for each $i=1,\ldots,m-1$ we have $\sum_{j\in D_i}y_j\geq (d+1)x_{i+1}$. This implies that
$$\sum_{i=1}^n y_i - \sum_{j\in D_m}y_j\ge (d+1) \sum_{i=2}^m x_i.$$
By Lemma~\ref{lem:o(1)}, $x_1$ and  $\sum_{j\in D_m}y_j$ are both $o(1)$.
It follows  that
$$\langle\x,\1\rangle=\sum_{i=1}^m x_i=o(1).$$
So by Lemma~\ref{lem:remarkdelta},
$$(1+o(1)) \frac{\pi^2}{m^2}=\mu(P_m) \le(1+o(1))\frac{ \sum_{i=1}^{m-1}(x_i-x_{i+1})^2}{\sum_{i=1}^m x_i^2}.$$
Now, by \eqref{eq:lowereven1} and \eqref{eq:lower}  we have
\begin{align*}
	\mu({\Gamma}) &=\frac{\sum_{ij\in E({\Gamma})}(y_i-y_j)^2}{\sum_{i=1}^n {y_i}^2 }\\
	&\ge(1+o(1))\frac\ell{(d+1)^2}\frac{ \sum_{i=1}^{m-1}(x_i-x_{i+1})^2}{\sum_{i=1}^m x_i^2}\\
	&= (1+o(1)) \frac\ell{(d+1)^2}\frac{\pi^2}{m^2}\\
	&=(1+o(1)) \frac{\ell \pi^2}{n^2}.
\end{align*}
\end{proof}

Now, we deduce that the upper and lower  bounds
given in Theorems~\ref{thm:upper} and \ref{thm:lower} can be extended to the graphs in ${\mathscr F}_{n,d,C}$.

\begin{theorem}\label{thm:upper-lower}
Let $\G\in {\mathscr F}_{n,d,C}$ with  major subgraphs 
 $\G(a_i,b_i,c_i;m_i)$, $i=1,\ldots,t$.
If $L$ and $\ell$ are the maximum and minimum of
 $\{a_ib_ic_i : i=1,\ldots,t\}$, respectively, then
$ (1+o(1))\frac{\ell\,\pi^2}{n^2} \le \mu(\G) \le (1+o(1))\frac{L\,\pi^2}{n^2}$.
\end{theorem}
\begin{proof}
 By the assumption  $\cal G$  is made of the major subgraphs $\G_i:=\G(a_i,b_i,c_i;m_i)$,    $i=1,\ldots,t$,  and some subgraphs $H_0,\ldots, H_t$  with $\sum_{i=0}^{t} |V(H_i)|\le C$. 
We let $$\g=\g_d(a_1,b_1,c_1,\ldots,a_t,b_t,c_t;m_1,\ldots,m_t).$$
By Theorems~\ref{thm:upper} and \ref{thm:lower}, we have $	 (1+o(1))\frac{\ell\pi^2}{n^2}\le\mu({\Gamma})\le(1+o(1))\frac{L\pi^2}{n^2}$. Note that $\G_1,\ldots,\G_t$ are also subgraphs of $\g$. We modify $\g$ to obtain $\G$ and show that this does not alter the order of the algebraic connectivity.
 
 We begin by incorporating the subgraphs $H_0,\ldots, H_t$ into $\g$, connecting vertices from the first cell $K_{a_i}$ and the last cell $K_{c_i}$ of $\G_i$ in $\g$ to $H_{i-1}$ and $H_i$, respectively, mirroring the edges between $H_{i-1},\G_i, H_i$ in $\G$. Let $\G'$ denote the resulting graph. Given that $t\le C$ and $H_i$, $K_{a_i}$, $K_{c_i}$ are all of order $O(1)$, applying Theorem~\ref{thm:H}, $t+1$ times, we conclude that
 $\mu(\G')=(1+o(1))\mu(\g)$.

Now to obtain $\G$ from $\G'$, we eliminate all edges between the last cell of $\G_i$ and the first cell of $\G_{i+1}$, for $i=1,\ldots,t-1$. It is evident that $\mu(\G)\leq \mu(\G')$. Hence, $\mu(\G)=o(1/n)$. Let $\x$ be a unit Fiedler vector of $\G$. 
Any pair of vertices adjacent in $\G'$ might not be adjacent in $\G$, but their distance in $\G$ is $O(1)$. Thus, by applying Lemma~\ref{lem:muedge}, for any $ij\in E(\G')\setminus E(\G)$, we have $(x_i-x_j)^2=o(\mu(\G))$.
Since $|E(\G')\setminus E(\G)|\le t d^2 \le C d^2$, it follows that $\x^\top L(\G')\x-\x^\top L(\G)\x=o(\mu(\G))$. This implies $\mu(\G')\leq\x^\top L(\G')\x=(1+o(1))\mu(\G)$. Therefore, we establish $\mu(\G)=(1+o(1))\mu(\G')$, and subsequently $\mu(\G)=(1+o(1))\mu(\g)$, from which the result follows.
\end{proof}

An immediate   consequence  of Theorem~\ref{thm:upper-lower} is the following corollary.

\begin{corollary}\label{coro:equal}
	Let $\G\in{\mathscr F}_{n,d,C}$ such that  its major  subgraphs are  all $\G(a,b,c;m)$. Then  $\mu(\G)=(1+o(1))\frac{a b c\,\pi^2}{n^2}$.
\end{corollary}

\section{Aldous--Fill and Guiduli--Mohar conjectures}\label{sec:proofs}

In this section, we present the proofs of Theorems~\ref{thm:ours=aldous} and \ref{thm:muDiamter3nd+1}, which we restate here for the reader's convenience.

 \firsttheorem*
\begin{proof}
  The graphs of Conjecture~\ref{conjecture:ours} belong to ${\mathscr F}_{n,d,C}$ with the major subgraph $\G(1,d-1,1;m)$ and $\G(1,2,d-2;m)$ for odd and even $d$, respectively.
So Corollary~\ref{coro:equal} implies that the algebraic connectivity  of these graphs   is equal to
$(1+o(1))\frac{(d-1)\pi^2}{n^2}$   for odd $d$  and $(1+o(1))\frac{2(d-2)\pi^2}{n^2}$ for  even $d$.
Therefore, if Conjecture~\ref{conjecture:ours} is  true, then for fixed $d\geq 3$, the maximum relaxation time over the family of $d$-regular graphs is $(1+o(1))\frac{dn^2}{(d-1)\pi^2}$ and $(1+o(1))\frac{dn^2}{2(d-2)\pi^2}$ for odd and even $d$, respectively.   Note that  for $x\geq 3$ the maximum value of the function $\frac{x}{x-1}$ is $\frac{3}{2}$ and for $x\geq4$ the maximum value of the function $\frac{x}{2(x-2)}$ is 1.
 So if  Conjecture~\ref{conjecture:ours} is  true, then the minimum algebraic connectivity and the maximum relaxation time over the family of all  regular graphs  with $n$ vertices is equal to $(1+o(1))\frac{2\pi^2}{n^2}$ and $(1+o(1)) \frac{3n^2}{2\pi^2}$, respectively; and are achieved by cubic graphs.
\end{proof}

\begin{figure}[t]
	\centering
	\begin{tikzpicture}[scale=.8]
		\draw (1,0) ellipse (.25 and 1.28);
		\vertex[fill] (1) at (-.3,0) [] {};
		\vertex[fill] (5) at (2.3,0) [] {};
		\vertex[fill] (2) at (1,1) [] {};
		\vertex[fill] (3) at (1,.6) [] {};
		\vertex[fill] (4) at (1,-1)[] {};
		\tikzstyle{vertex}=[circle, draw, inner sep=.3pt, minimum size=.3pt]
		\vertex[fill] () at (1,-.1) [] {};
		\vertex[fill] () at (1,-.2) [] {};
		\vertex[fill] () at (1,-.3) [] {};
		\tikzstyle{vertex}=[circle, draw, inner sep=0pt, minimum size=0pt]
		\vertex[] () at (1,1.3) [label=above:\footnotesize{$K_{d-1}$}] {};
		\vertex[] (0) at (-.6,0) [] {};
		\vertex[] (6) at (2.6,0) [] {};
		\path
		(1) edge (2)
		(1) edge (3)
		(1) edge (4)
		(5) edge (2)
		(5) edge (3)
		(1) edge (0)
		(5) edge (6)
		(5) edge (4);
	\end{tikzpicture}	\quad \quad \quad
	\begin{tikzpicture}[scale=.8]
		\draw (-1,0) ellipse (.25 and 1.28);
		\vertex[fill] (1) at (.1,0) [] {};
		\vertex[fill] (5) at (-2.3,.4) [] {};
		\vertex[fill] (55) at (-2.3,-.4) [] {};
		\vertex[fill] (555) at (-2.8,0) [] {};
		\vertex[fill] (2) at (-1,1) [] {};
		\vertex[fill] (3) at (-1,.6) [] {};
		\vertex[fill] (4) at (-1,-1) [] {};
		\tikzstyle{vertex}=[circle, draw, inner sep=0pt, minimum size=0pt]
		\vertex[] () at (-1,1.3) [label=above:\footnotesize{$K_{d-2}$}] {};
		\vertex[fill] (7) at (-3.1,.27) [ ] {};  
		\vertex[fill] (7a) at (-3.1,.15) [ ] {}; 
		\vertex[fill] (77) at (-3.1,-.27) [ ] {};
		\vertex[fill] (777) at (.4,.27) [ ] {};  
		\vertex[fill] (7777) at (.4,-.27) [ ] {};
		\tikzstyle{vertex}=[circle, draw, inner sep=.3pt, minimum size=.3pt]
		\vertex[fill] () at (-1,-.1) [] {};
		\vertex[fill] () at (-1,-.2) [] {};
		\vertex[fill] () at (-1,-.3) [] {};
		\tikzstyle{vertex}=[circle, draw, inner sep=.1pt, minimum size=.1pt]
		\vertex[fill] () at (-3.1,-.01) [ ] {};  
		\vertex[fill] () at (-3.1,-.06) [ ] {}; 
		\vertex[fill] () at (-3.1,-.11) [ ] {};
		\path
		(1) edge (777)
		(1) edge (7777) 
		(555) edge (7)
		(555) edge (77)
		(555) edge (7a)
		(1) edge (2)
		(1) edge (3)
		(1) edge (4)
		(5) edge (2)
		(5) edge (3)
		(5) edge (4)
		(55) edge (2)
		(55) edge (3)
		(55) edge (4)
		(55) edge (5)
		(555) edge (5)
		(555) edge (55)	;
	\end{tikzpicture}\caption{The blocks $L_d$ (left) and $M_d$ (right).}\label{fig:Type}
\end{figure}
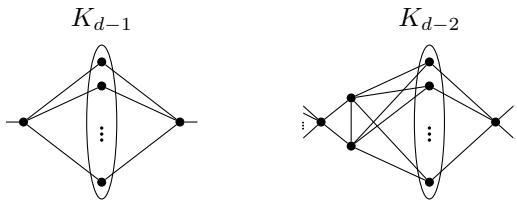

Now we prove Theorem~\ref{thm:muDiamter3nd+1}
as a consequence of Theorem~\ref{thm:upper-lower}.
 As shown in Figure~\ref{fig:Type}, we denote the blocks $K_1+K_{d-1}+K_1$ and $K_1+K_2+K_{d-2}+K_1$ by $L_d$ and $M_d$, respectively.
 \secondtheorem*
\begin{proof}
By Theorem~\ref{thm:diam}, it is enough to prove the assertion for 
 the graphs in  ${\mathscr F}_{n,d,C}$. 
 
(i) Let $\G\in{\mathscr F}_{n,d,C}$ and $\G_1,\ldots, \G_t$ be  the  major subgraphs of $\G$. By Theorem~\ref{thm:upper-lower},
$ (1+o(1))\frac{\ell\,\pi^2}{n^2}\le \mu(\G)$, where $\ell:=\min \{a_ib_ic_i : i=1,\ldots,t\}$. The minimum of the function $f(x,y,z)=xyz$, subject to $x+y+z=d+1$ and $x,y,z\geq 1$ is $d-1$. This means that
$\ell=d-1$. On the other hand, by Corollary~\ref{coro:equal},
the path-like graph  with $m=\lfloor n/(d+1)\rfloor$ blocks $L_d$ (of Figure~\ref{fig:Type}) attains the minimum $\mu=(1+o(1))\frac{(d-1)\,\pi^2}{n^2}$.

(ii)
Let $\G$ be a $d$-regular graph with minimum $\mu$ in  ${\mathscr F}_{n,d,C}$. 
By (i), for odd $d$,  $\mu(\G)=(1+o(1))\frac{(d-1)\pi^2}{n^2}$.
Let $d$  be  even and  $\G_1,\ldots, \G_t$ be  the  major subgraphs of $\G$. By Theorem~\ref{thm:upper-lower},
$ (1+o(1))\frac{\ell\,\pi^2}{n^2}\le \mu(\G)$, where $\ell:=\min \{a_ib_ic_i : i=1,\ldots,t\}$.
 For positive integers $x,y,z$, the  minimum  of  $f(x,y,z)=xyz$ subject to  $x+y+z=d+1$,  and $x+y,x+z,y+z\ge3$ (this condition is necessary as $\G$ has no bridge) occurs   if  $x,y,z$ are $1,2,d-2$ in any order. 
Thus $\ell=2(d-2)$.
On the other hand, a path-like $d$-regular graph whose  blocks (except the end ones) are $M_d$ (of Figure~\ref{fig:Type})
 attains the minimum $\mu=(1+o(1))\frac{2(d-2)\pi^2}{n^2}$.
 For odd (resp., even) values of $d$, examples of $d$-regular graphs with all middle blocks $L_d$ (resp., $M_d$) are provided in Table~\ref{TableMax}.

The rest of the assertion follows immediately.
\end{proof}

\section{Max diameter versus min algebraic connectivity}\label{sec:MaxMin}

In our last section, we investigate the interplay between the graphs with maximum diameter and those with minimum algebraic connectivity within the family of $d$-regular graphs or those with $\delta=d$. 
In this regard, we find it natural to consider the following extension of Problem~\ref{problem1}:
\begin{itemize}
	\item[(a)] Given $d\ge3$, is it true that among $d$-regular graphs (or graphs with $\delta=d$) those with maximum diameter have algebraic connectivity smaller than  others?	
\end{itemize}
We shall see that the answer to this question is negative. So one might wonder whether its asymptotic variation holds:
\begin{itemize}
	\item[(b)] Given $d\ge3$, is it true that among $d$-regular graphs (or graphs with $\delta=d$) those with  asymptotically maximum diameter have  asymptotically minimum algebraic connectivity?
\end{itemize}
 Based on the results of Section~\ref{sec:MaxDiam}, we address this variation as well, with the exception of $3$- and $4$-regular graphs and graphs with $\delta=3$. For $3$- and $4$-regular graphs, we present a weaker variant applicable to those with diameter  $\frac{3n}{d+1}+O(1)$. We propose the converse of (b) as a conjecture.

The diameter of a graph can be bounded in terms of its order and minimum degree.
Several results in this line can be found in the literature  (see, e.g., \cite{Caccetta,Erdos,Moon,Soares}). The first result of this type can be attributed to Moon \cite{Moon},
who proved that for a  graph $G$ of order $n$ and minimum degree $d \geq 2$, $\dm(G)\leq (3n-2d-6)/d$.
The following result determines the maximum diameter explicitly. 
\begin{theorem}[Caccetta and  Smyth \cite{Caccetta}]\label{thm:diameter}
 The maximum diameter of a  graph of order $n$ and minimum degree $d$
	\begin{itemize}
		\item[\rm(i)] for $n\leq 2d+1$ is $\left\lceil\frac{n}{d+1}\right\rceil$,
		\item[\rm(ii)] for $n \geq 2d+2$  is  $3\left\lfloor\frac{n}{d+1}\right\rfloor-\left\{\begin{array}{l}
			3 \quad n\equiv0\pmod{d+1},\\
			2 \quad n\equiv1\pmod{d+1},\\
			1 \quad \text{otherwise.}\end{array} \right.$
	\end{itemize}
\end{theorem}
In \cite{Caccetta}, it was also shown  that no $d$-regular graph of diameter $3\left\lfloor\frac{n}{d+1}\right\rfloor-1$  exists when $d$ is even and $n\equiv2\pmod{d+1}$. In this case,  we observe that $d$-regular graphs of diameter $3\left\lfloor\frac{n}{d+1}\right\rfloor-2$ exist.  Thus,  the  following theorem can be deduced.
 \begin{theorem}\label{thm:diameterReg}
	%For $n\geq2d+4$ and $G\in {\mathfrak{D}}_{n,d}$, the following statements are valid:
Let $d\geq3$ and  $n\geq2d+4$. The maximum diameter of a $d$-regular graph of order $n$
	\begin{itemize}
		\item[\rm(i)] for odd $d$  is $3\left\lfloor\frac{n}{d+1}\right\rfloor-\left\{\begin{array}{l}
			3 \quad n\equiv0\pmod{d+1},\\
			1 \quad \text{otherwise,}\end{array} \right.$
		\item[\rm(ii)] for even $d$ is  $3\left\lfloor\frac{n}{d+1}\right\rfloor-\left\{\begin{array}{l}
			3 \quad n\equiv0\pmod{d+1},\\
			2 \quad n\equiv1,2\pmod{d+1},\\
			1 \quad \text{otherwise.}\end{array} \right.$
	\end{itemize}
\end{theorem}

\begin{table}
	\centering
	\scalebox{.62}{
		\begin{tabular}{|c|c|c|l}\cline{1-3}
			$d$ & $r$ & All middle blocks are $L_d$ or $M_d$  \\ \hline
			& & \\
			odd & 0  & $K_2+K_{d-1}^{-1}+(K_1+K_1+K_{d-1})_{m-3}+K_1+K_1+K_{d-1}\overset{{\scriptscriptstyle +1}}{\cup} K_{d-1}^{-1}+K_2$\\
			& $2,4,\ldots,d-1$ &  $K_2+K_{d-1}^{-1}+(K_1+K_1+K_{d-1})_{m-2}+K_1+K_1+K_{d-1}^{-(r-1)}+ \overline{ K}_{r}^{+1}$\\
			& &  \\
			& 0 &   $K_3+K_{d-2}^{-1}+(K_1+K_2+K_{d-2})_{m-3}+K_1+K_2+K_{d-2}\overset{{\scriptscriptstyle +1}}{\cup} K_{d-2}^{-1}+K_3$\\
			& 1 &  $K_3+K_{d-2}^{-1}+(K_1+K_2+K_{d-2})_{m-2}+K_1+\overline{K}_2+K_{d-1}$\\
			even  & 2 & $K_3+K_{d-2}^{-1}+(K_1+K_2+K_{d-2})_{m-2}\circ \mathcal{H}_1$\\
			& $3,4,\ldots,d-2$ &  $K_3+K_{d-2}^{-1}+(K_1+K_2+K_{d-2})_{m-2}+K_1+K_2+ K_{d-2}^{-(r-1)}+C_r$\\
			& $d-1$ &    $K_3+K_{d-2}^{-1}+(K_1+K_2+K_{d-2})_{m-2}+K_1+\overline{K}_2+\overline{K}_{d-1}+ \overline{K}_{d-2}^{+1}$\\
			& $d$ &   $K_3+K_{d-2}^{-1}+(K_1+K_2+K_{d-2})_{m-2}+K_1+\overline{K}_2+\overline{K}_{d-1}\overset{-1}{+}C_{d-1}$\\
			& &  \\
			\hline
	\end{tabular}}\caption{Some members of $\mathscr{D}_{n,d}$; here
		$n=(d+1)m+r$ with $m\ge3$ and $0\le r\le d$.}
	\label{TableMax}
\end{table}

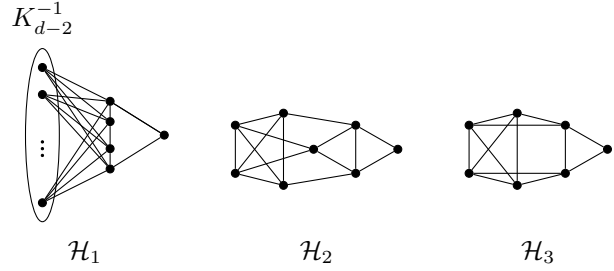
\begin{figure}
	\captionsetup[subfigure]{labelformat=empty}
	\centering
	\subfloat[$\mathcal{H}_1$]{\begin{tikzpicture}[scale=.9]
			\draw (-1.5,0) ellipse (.25 and 1.28);
			\vertex[fill] (1) at (.3,0) [] {};
			\vertex[fill] (2) at (-.5,.5) [] {};
			\vertex[fill] (3) at (-.5,.2) [] {};
			\vertex[fill] (4) at (-.5,-.2) [] {};
			\vertex[fill] (5) at (-.5,-.5) [] {};
			\vertex[fill] (6) at (-1.5,1) [] {};
			\vertex[fill] (7) at (-1.5,.6) [] {};
			\vertex[fill] (8) at (-1.5,-1) [] {};
			\tikzstyle{vertex}=[circle, draw, inner sep=.3pt, minimum size=.3pt]
			\vertex[fill] () at (-1.5,-.1) [] {};
			\vertex[fill] () at (-1.5,-.2) [] {};
			\vertex[fill] () at (-1.5,-.3) [] {};
			\tikzstyle{vertex}=[circle, draw, inner sep=0pt, minimum size=0pt]
			\vertex[] () at (-1.5,1.3) [label=above:\footnotesize{$K^{-1}_{d-2}$}] {};
			\path
			(1) edge (2)
			(1) edge (5)
			(1) edge (2)
			(2) edge (3)
			(3) edge (4)
			(4) edge (5)
			(2) edge (6)
			(2) edge (7)
			(2) edge (8)
			(3) edge (6)
			(3) edge (7)
			(3) edge (8)
			(4) edge (6)
			(4) edge (7)
			(4) edge (8)
			(5) edge (6)
			(5) edge (7)
			(5) edge (8) ;
	\end{tikzpicture}}
	\qquad
	\subfloat[$\mathcal{H}_2$]{\begin{tikzpicture}[scale=0.8]
			\vertex[fill] (r1) at (0,0) [] {};
			\vertex[fill] (r2) at (0,.8) [] {};
			\vertex[fill] (r3) at (.8,1) [] {};
			\vertex[fill] (r4) at (.8,-.2) [] {};
			\vertex[fill] (r5) at (1.3,.4) [] {};
			\vertex[fill] (r6) at (2,.8) [] {};
			\vertex[fill] (r7) at (2,0) [] {};
			\vertex[fill] (r8) at (2.7,.4) [] {};
			\tikzstyle{vertex}=[circle, draw, inner sep=0pt, minimum size=0pt]
			\vertex[] () at (0,-.8) [ ] {};
			\path
			(r5) edge (r2)
			(r1) edge (r2)
			(r1) edge (r5)
			(r1) edge (r3)
			(r1) edge (r4)
			(r2) edge (r3)
			(r2) edge (r4)
			(r3) edge (r6)
			(r4) edge (r3)
			(r4) edge (r7)
			(r5) edge (r7)
			(r5) edge (r6)
			(r6) edge (r7)
			(r7) edge (r8)
			(r6) edge (r8)  ;
	\end{tikzpicture}}
	\qquad
	\subfloat[$\mathcal{H}_3$]{\begin{tikzpicture}[scale=.8]
			\vertex[fill] (r1) at (0,.8) [] {};
			\vertex[fill] (r2) at (.8,1) [] {};
			\vertex[fill] (r3) at (0,0) [] {};
			\vertex[fill] (r4) at (.8,-.2) [] {};
			\vertex[fill] (r5) at (1.6,.8) [] {};
			\vertex[fill] (r6) at (1.6,0) [] {};
			\vertex[fill] (r7) at (2.3,.4) [] {};
			\tikzstyle{vertex}=[circle, draw, inner sep=0pt, minimum size=0pt]
			\vertex[] () at (0,-.8) [ ] {};
			\path
			(r5) edge (r2)
			(r1) edge (r2)
			(r1) edge (r5)
			(r1) edge (r3)
			(r1) edge (r4)
			(r2) edge (r3)
			(r2) edge (r4)
			(r3) edge (r6)
			(r4) edge (r3)
			(r4) edge (r6)
			(r5) edge (r6)
			(r7) edge (r6)
			(r5) edge (r7);
	\end{tikzpicture}}
	\caption{Three  possible  end blocks for  regular graphs of maximum diameter.}
	\label{fig:HMax}
\end{figure}

We denote  the family of $d$-regular graphs with $n$ vertices and maximum diameter by $\mathscr{D}_{n,d}$.
The graphs in $\mathscr{D}_{n,3}$ has been characterized in \cite{Imrich}: path-like graphs all whose middle blocks are $L_3$ (see Figure~\ref{fig:Type}).
From Theorem~\ref{thm:diam} and its proof (also from \cite{AbGh}),
it is not hard to understand the structure of  the graphs in $\mathscr{D}_{n,4}$.
In fact, such a graph  has a path-like structure and almost every three consecutive parts in its distance partition  together  have $5$ vertices. Then the regularity condition implies that all blocks (with few exceptions) are $M_4$.
For general $d$, some  members of $\mathscr{D}_{n,d}$  are identified in Table~\ref{TableMax}.
The notation used in this table is clarified below.
%Here, we find it more intuitive to use the notation $(K_a+K_b+K_c)_m$ to represent the same graph as %$\G(a,b,c;m)$.
 %{and by \R $(K_a+K_b+K_c)_m+K_d$ or $K_d+(K_a+K_b+K_c)_m$,
%we mean the graph obtained by joining each vertex of $K_d$
%with each vertex of  the last clique  $K_c$ or  the first clique $K_a$, respectively, of $(K_a+K_b+K_c)_m$.}
As usual,  $C_n$  denotes the  {\em cycle} of length $n$ and $\overline{G}$  the {\em complement} of  $G$.
%{\R and $G+H$ the join of two graphs $G$ and $H$.}
%A {\em $1$-factor} of a graph $G$ is a $1$-regular spanning subgraph of $G$.
%$K_{2m}$ has $2m-1$ $1$-factors.
By $G^{-r}$ (resp., $G^{+r}$) we mean the graph obtained from $G$ by removing (resp., adding)  the edges of
$r$ $1$-factors.
When $|V(G)|=|V(H)|$, denote by $G\overset{{\scriptscriptstyle +1}}{\cup}H$ (resp., $G\overset{{\scriptscriptstyle  -1}}{+}H$) the graph obtained from $G\cup H$ (resp., $G+H$)  by adding (resp., removing) edges of one $1$-factor between  $G$ and $H$. 
In a sequential join of graphs, when some of the summands are repeated, for example, in the case of $G+K_a+K_b+K_c+\cdots+K_a+K_b+K_c+H$, where $K_a+K_b+K_c$ is repeated $m$ times, we use the notation $G+(K_a+K_b+K_c)_m+H$ for brevity.
Finally, given the graphs $\mathcal{H}_i$ shown in Figure~\ref{fig:HMax}, by  $(K_a+K_b+K_c)_m\circ \mathcal{H}_i$ or $\mathcal{H}_i\circ(K_a+K_b+K_c)_m$,
we mean the graph obtained by joining the vertex of degree 2 in $\mathcal{H}_i$ 
to the last clique  $K_c$ or  the first clique  $K_a$, respectively, of $(K_a+K_b+K_c)_m$. 

%{\R Let $i\in \{1,\ldots, t\}$ and $\ast_i\in\{+,\overset{{\scriptscriptstyle +1}}{\cup},\overset{{\scriptscriptstyle  -1}}{+}\}$, in the sequence $G_1\ast_1G_2\ast_2\cdots\ast_t G_{t+1}$ each $\ast_i$ operates only between the graphs $G_i$ and $G_{i+1}$. If $G=G_1\ast_1\cdots\ast_t G_{t+1}$ to simplify instead of the squence $G+\cdots+G$ consists of $m$ graphs $G$, we write $(G)_m$. 

Now, we are prepared to prove the final theorem of the paper. Part~(i) provides a negative answer to (a), particularly addressing Problem~\ref{problem1}. Part~(iii) demonstrates that (b) fails for $d$-regular graphs, as well as graphs with $\delta \ge d$ for $d \ge 5$, and Part~(iv) establishes the same for graphs with $\delta = 4$. The correctness of (b) for $3$- and $4$-regular graphs, and graphs with $\delta = 3$, remains an open question. Though Part~(ii) establishes a weaker version for nearly maximum-diameter graphs.

\begin{theorem}\label{thm:diam-vs-mu}
	\begin{itemize}
		\item[\rm(i)] For every $d\ge3$, for some $n$, there  exist $n$-vertex  $d$-regular graphs $\g$ and $\g'$ such that $\g\in\mathscr{D}_{n,d}$ and $\g'\not\in\mathscr{D}_{n,d}$  but $\mu(\g')<\mu(\g)$.
		\item[\rm(ii)] For $d=3,4$,  $d$-regular graphs with  diameter $\frac{3n}{d+1}+O(1)$ have asymptotically minimum algebraic connectivity.
		\item[\rm(iii)] For any $d\ge5$, there are sequences of $d$-regular graphs $\g_n$ of asymptotically  maximum diameter and $\g'_n$ with $\dm(\g'_n)<(1-\e)\dm(\g_n)$ such that $\mu(\g'_n)<(1-\e)\mu(\g_n)$ for some $\e>0$.
		\item[\rm(iv)] There are graphs with $\delta=4$ and  asymptotically maximum diameter that do  not have asymptotically minimum algebraic connectivity.
	\end{itemize}
\end{theorem}

\begin{proof}
	(i) Let $m$ be even, $n=4m+16$, and
	$$\g_n=K_2+K_{2}^{-1}+(K_1+K_1+K_{2})_{\frac{m}{2}}+K_1+K_1+K_{2}\overset{{\scriptscriptstyle +1}}{\cup} K_{2}+
	(K_1+K_1+K_{2})_{\frac{m}{2}}+K_1+K_1+ K_{2}^{-1}+K_2,$$
	$$\g'_n=K_2+K_{2}^{-1}+(K_1+K_1+K_{2})_{m}+K_1+K_1+K_{2}\overset{{\scriptscriptstyle +1}}{\cup}K_{2}\overset{{\scriptscriptstyle +1}}{\cup}K_{2}\overset{{\scriptscriptstyle +1}}{\cup} K_{2}^{-1}+K_2.$$
 See Figure~\ref{fig:mincubic} for an illustration of these two graphs. We observe that $\dm(\g_n)=3m+9=3n/4-3$ and thus by Theorem~\ref{thm:diameterReg}, $\g_n\in \mathscr{D}_{n,3}$. Also  $\dm(\g'_n)=\dm(\g_n)-1$.
	Using  computer, we observed that for quit a few values of $n$, for instance any $n=4m+16$ with $4\le m\le260$, we have   $\mu(\g_n')<\mu(\g_n)$.\footnote{We believe that this is true for every $m \geq 4$. A rigorous proof involves tedious calculations, which we do not pursue here.} 
	Similarly for quartic graphs, let $n=5m+13$ and
	$$\g_n=K_3+K_{2}^{-1}+(K_1+K_2+K_{2})_{m}\circ \mathcal{H}_2,$$
$$\g'_n=\mathcal{H}_3 \circ (K_2+K_2+K_1)+K_2+K_2+ (K_1+K_2+K_{2})_{m-2}\circ \mathcal{H}_3,$$
	where  $\mathcal{H}_2$ and $\mathcal{H}_3$ are the graphs depicted in Figure~\ref{fig:HMax}.
	It is easy to verify that $\g\in \mathscr{D}_{n,4}$  and  $\dm(\g'_n)=\dm(\g_n)-1$. 
		Again using  computer, we observed that for any $n=5m+13$ with  $1\le m\le 260$, we have   $\mu(\g_n')<\mu(\g_n)$.
	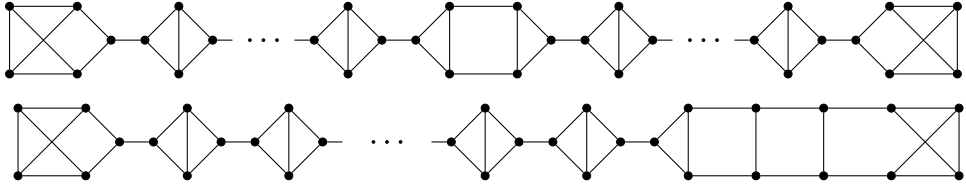
\begin{figure}
	\centering
	\begin{tikzpicture}[scale=0.9]
		\vertex[fill] (1) at (-1,-.5) [] {};
		\vertex[fill] (2) at (-1,.5) [] {};
		\vertex[fill] (3) at (0,-.5) [] {};
		\vertex[fill] (4) at (0,.5) [] {};
		\vertex[fill] (5) at (0.5,0) [] {};
		\vertex[fill] (6) at (1,0) [] {};
		\vertex[fill] (7) at (1.5,.5) [] {};
		\vertex[fill] (8) at (1.5,-.5) [] {};
		\vertex[fill] (9) at (2,0) [] {};
		
		\vertex[fill] (10) at (3.5,0) [] {};
		\vertex[fill] (11) at (4,.5) [] {};
		\vertex[fill] (12) at (4,-.5) [] {};
		\vertex[fill] (13) at (4.5,0) [] {};
		
		\vertex[fill] (a) at (5,0) [] {};
		\vertex[fill] (b) at (5.5,.5) [] {};
		\vertex[fill] (c) at (5.5,-.5) [] {};
		\vertex[fill] (d) at (6.5,.5) [] {};
		\vertex[fill] (e) at (6.5,-.5) [] {};
		\vertex[fill] (f) at (7,0) [] {};
		
		\vertex[fill] (22) at (7.5,0) [] {};
		\vertex[fill] (23) at (8,.5) [] {};
		\vertex[fill] (24) at (8,-.5) [] {};
		\vertex[fill] (25) at (8.5,0) [] {};
		
		\vertex[fill] (26) at (10,0) [] {};
		\vertex[fill] (27) at (10.5,.5) [] {};
		\vertex[fill] (28) at (10.5,-.5) [] {};
		\vertex[fill] (29) at (11,0) [] {};
		\vertex[fill] (34) at (11.5,0) [] {};
		\vertex[fill] (33) at (12,.5) [] {};
		\vertex[fill] (32) at (12,-.5) [] {};
		\vertex[fill] (31) at (13,.5) [] {};
		\vertex[fill] (30) at (13,-.5) [] {};
		\tikzstyle{vertex}=[circle, draw, inner sep=0pt, minimum size=1pt]
		\vertex[fill] () at (2.75,0) [] {};
		\vertex[fill] () at (2.95,0) [] {};
		\vertex[fill] () at (2.55,0) [] {};
		\vertex[fill] () at (9.25,0) [] {};
		\vertex[fill] () at (9.45,0) [] {};
		\vertex[fill] () at (9.05,0) [] {};
		\tikzstyle{vertex}=[circle, draw, inner sep=0pt, minimum size=0pt]
		\vertex[] (14) at (2.3,0) [] {};
		\vertex[] (21) at (3.2,0) [] {};
		\vertex[] (cc) at (8.8,0) [] {};
		\vertex[] (ccc) at (9.7,0) [] {};
		\path
		(1) edge (2)
		(1) edge (3)
		(1) edge (4)
		(2) edge (3)
		(2) edge (4)
		(3) edge (5)
		(4) edge (5)
		(5) edge (6)
		(6) edge (7)
		(6) edge (8)
		(7) edge (8)
		(7) edge (9)
		(8) edge (9)
		(9) edge (14)
		(10) edge (21)
		(10) edge (11)
		(10) edge (12)
		(11) edge (12)
		(11) edge (13)
		(12) edge (13)
		(13) edge (a)
		(a) edge (b)
		(a) edge (c)
		(b) edge (c)
		(b) edge (d)
		(e) edge (c)
		(d) edge (e)
		(d) edge (f)
		(e) edge (f)
		(cc) edge (25)
		(ccc) edge (26)
		(f) edge (22)
		(22) edge (23)
		(22) edge (24)
		(23) edge (24)
		(23) edge (25)
		(24) edge (25)
		%%%%(25) edge (26)
		(26) edge (27)
		(26) edge (28)
		(27) edge (28)
		(27) edge (29)
		(28) edge (29)
		(29) edge (34)
		(30) edge (31)
		(30) edge (32)
		(31) edge (33)
		(32) edge (31)
		(30) edge (33)
		(32) edge (34)
		(33) edge (34);
	\end{tikzpicture}
	\vspace{.3cm} \\
	\begin{tikzpicture}[scale=0.9]
		\vertex[fill] (1) at (0,-.5) [] {};
		\vertex[fill] (2) at (0,.5) [] {};
		\vertex[fill] (3) at (1,-.5) [] {};
		\vertex[fill] (4) at (1,.5) [] {};
		\vertex[fill] (5) at (1.5,0) [] {};
		\vertex[fill] (6) at (2,0) [] {};
		\vertex[fill] (7) at (2.5,.5) [] {};
		\vertex[fill] (8) at (2.5,-.5) [] {};
		\vertex[fill] (9) at (3,0) [] {};
		\vertex[fill] (10) at (3.5,0) [] {};
		\vertex[fill] (11) at (4,.5) [] {};
		\vertex[fill] (12) at (4,-.5) [] {};
		\vertex[fill] (13) at (4.5,0) [] {};
		\vertex[fill] (22) at (6.4,0) [] {};
		\vertex[fill] (23) at (6.9,.5) [] {};
		\vertex[fill] (24) at (6.9,-.5) [] {};
		\vertex[fill] (25) at (7.4,0) [] {};
		\vertex[fill] (26) at (7.9,0) [] {};
		\vertex[fill] (27) at (8.4,.5) [] {};
		\vertex[fill] (28) at (8.4,-.5) [] {};
		\vertex[fill] (29) at (8.9,0) [] {};
		\vertex[fill] (34) at (10.9,-.5) [] {};
		\vertex[fill] (33) at (10.9,.5) [] {};
		\vertex[fill] (32) at (9.9,-.5) [] {};
		\vertex[fill] (31) at (9.9,.5) [] {};
		\vertex[fill] (30) at (9.4,0) [] {};
		\vertex[fill] (35) at (11.9,-.5) [] {};
		\vertex[fill] (36) at (11.9,.5) [] {};
		\vertex[fill] (37) at (12.9,-.5) [] {};
		\vertex[fill] (38) at (12.9,.5) [] {};
		\vertex[fill] (39) at (13.9,-.5) [] {};
		\vertex[fill] (40) at (13.9,.5) [] {};
		
		\tikzstyle{vertex}=[circle, draw, inner sep=0pt, minimum size=1pt]
		\vertex[fill] (15) at (5.25,0) [] {};
		%\vertex[fill] (16) at (5.3,0) [] {};
		\vertex[fill] (17) at (5.45,0) [] {};
		%\vertex[fill] (18) at (5.5,0) [] {};
		\vertex[fill] (19) at (5.65,0) [] {};
		%\vertex[fill] (20) at (5.7,0) [] {};
		\tikzstyle{vertex}=[circle, draw, inner sep=0pt, minimum size=0pt]
		\vertex[] (14) at (4.8,0) [] {};
		\vertex[] (21) at (6.1,0) [] {};
		\path
		(1) edge (2)
		(1) edge (3)
		(1) edge (4)
		(2) edge (3)
		(2) edge (4)
		(3) edge (5)
		(4) edge (5)
		(5) edge (6)
		(6) edge (7)
		(6) edge (8)
		(7) edge (8)
		(7) edge (9)
		(8) edge (9)
		(9) edge (10)
		(10) edge (11)
		(10) edge (12)
		(11) edge (12)
		(11) edge (13)
		(12) edge (13)
		(13) edge (14)
		%(15) edge (16)
		%(17) edge (18)
		%(19) edge (20)
		(21) edge (22)
		(22) edge (23)
		(22) edge (24)
		(23) edge (24)
		(23) edge (25)
		(24) edge (25)
		(25) edge (26)
		(26) edge (27)
		(26) edge (28)
		(27) edge (28)
		(27) edge (29)
		(28) edge (29)
		(29) edge (30)
		(30) edge (31)
		(30) edge (32)
		(31) edge (32)
		(31) edge (33)
		(32) edge (34)
		(33) edge (34)
		(33) edge (36)
		(34) edge (35)
		(36) edge (35)
		(38) edge (36)
		(37) edge (35)
		(39) edge (37)
		(37) edge (40)
		(38) edge (39)
		(38) edge (40)
		(39) edge (40)
		;
	\end{tikzpicture}
	\caption{The cubic graphs $\Gamma$ (top) and $\Gamma'$ (bottom) of the proof of Theorem~\ref{thm:diam-vs-mu}\,(i).}
	\label{fig:mincubic}
\end{figure}	

 Now, suppose that  $d\geq 5$ be odd, $m\geq (d+7)/2$, $n=m(d+1)$,  and $$\g_n=K_{4}+K_{d-3}^{-1}+(K_1+K_3+K_{d-3})_{m-3}+K_1+K_3+K_{d-3}\overset{{\scriptscriptstyle  +1}}{\cup} K_{d-3}^{-1}+K_{4}.$$
		We have $\dm(\g_n)=3m-3=\frac{3n}{d+1}-3$ and thus by Theorem~\ref{thm:diameterReg}, $\g_n\in \mathscr{D}_{n,d}$. Furthermore,   by Corollary~\ref{coro:equal},  $\mu(\Gamma_n)=(1+o(1))\frac{3(d-3)\pi^2}{n^2}$. 
Consider the following graph, also from $\mathscr{D}_{n,d}$:
		$$G_n=K_2+K_{d-1}^{-1}+(K_1+K_1+K_{d-1})_{m-3}+K_1+K_1+K_{d-1}\overset{{\scriptscriptstyle +1}}{\cup} K_{d-1}^{-1}+K_2.$$
%We see that $\dm(G_n)=3m-3=\dm(\Gamma_n)$.		
		In $G_n$, replace a subgraph
		$(K_1+K_1+K_{d-1})_{\frac{d+1}{2}}$ by the subgraph
$$K_1+K_1+K_{d-1}\overset{{\scriptscriptstyle +1}}{\cup}K_{d-1}\overset{{\scriptscriptstyle +1}}{\cup}\cdots \overset{{\scriptscriptstyle +1}}{\cup}K_{d-1},$$
 consists of $(d+7)/2$ cells.
		Thus for the resulting graph $\Gamma'_n$, we have  $\dm(\Gamma'_n)=3m-d-1$.
	By   Corollary~\ref{coro:equal}, $\mu(G_n)=(1+o(1))\frac{(d-1)\pi^2}{n^2}$.
		As $\g'_n$ is obtained from $G_n$ by an $O(1)$-perturbation, from
		Theorem~\ref{thm:H} it follows that		
		$\mu(\g'_n)=(1+o(1))\mu(G_n)=(1+o(1))\frac{(d-1)\pi^2}{n^2}$, and thus $\mu(\g'_n)$ is asymptotically smaller than $\mu(\g_n)$.

		Finally, suppose that  $d\geq 6$ be even, $m\geq 2(d+1)$, $n=m(d+1)+4$,  and 
		$$\g_n=K_{3}+K_{d-2}^{-2}+{\overline K}_{2}+K_2+K_{d-3}+(K_2+K_2+K_{d-3})_{m-3}+K_2+{\overline K}_{2}+ K_{d-2}^{-2}+K_{3}.$$
		We have $\dm(\g_n)=3m-1=3\lfloor\frac{n}{d+1}\rfloor-1$ and thus by Theorem~\ref{thm:diameterReg}, $\g_n\in \mathscr{D}_{n,d}$. Furthermore,   by Corollary~\ref{coro:equal},  $\mu(\Gamma_n)=(1+o(1))\frac{4(d-3)\pi^2}{n^2}$. 
Consider the following graph, also from $\mathscr{D}_{n,d}$:
		$$ G_n=K_3+K_{d-2}^{-1}+(K_1+K_2+K_{d-2})_{m-2}+K_1+K_2+ K_{d-2}^{-3}+C_4.$$
		In $G_n$, replace the subgraph
		$(K_1+K_2+K_{d-2})_{2d}+K_1+K_2$
		by the subgraph
$$K_1+K_2+(K_{d-2}\overset{{\scriptscriptstyle +1}}{\cup}K_{d-2}+{\overline K}_{2}+{\overline K}_{2})_{d+1}.$$
	For the resulting graph $\Gamma'_n$, we have  $\dm(\Gamma'_n)=3m-2d+3$.
 From	Theorem~\ref{thm:H} and Corollary~\ref{coro:equal}, it follows that		
$\mu(\g'_n)=(1+o(1))\mu(G_n)=(1+o(1))\frac{2(d-2)\pi^2}{n^2}$.
		So $\mu(\Gamma'_n)$ is asymptotically smaller than
		$\mu(\Gamma_n)$.

	(ii)  First consider $d=3$. 
	Let $G$ be a cubic graph with $\dm(G)=3n/4+O(1)$.
By Theorem~\ref{thm:diam}, for some constant $C$, $G$ belongs to the family $\mathscr{F}_{n,3,C}$, with major subgraphs $\G(a,b,c;m)$ where
$a+b+c=4$. However, the only possible solution for this equation is $1,1,2$ in any order. It follows that (cf. the proof of Theorem~\ref{thm:diam}) that all the middle blocks of $G$ with few exceptions must be $L_3$, and thus by Corollary~\ref{coro:equal}, $\mu(G)=(1+o(1))\frac{2\pi^2}{n^2}$. By \cite{AbGhIm},  this is in fact minimum $\mu$ of cubic graphs. 

	Next, assume that $G$ is a quartic graph with $\dm(G)=3n/5+O(1)$.
	For $d=4$, we should find the solutions of  $a+b+c=5$,  subject to $a+b,a+c,b+c\ge3$ (since $G$ should have no bridge). It follows that  $a,b,c$ are $1,2,2$ in any order. So the middle blocks of $G$ with few exceptions must be $M_4$ and thus by Corollary~\ref{coro:equal}, $\mu(G)=(1+o(1))\frac{4\pi^2}{n^2}$.
	By \cite{AbGh}, this is minimum $\mu$ of quartic graphs.
	
	(iii)  Let $d\ge5$ and $\g_n$  be a $d$-regular path-like graph all whose middle blocks are $K_2+K_{d-3}+K_2$. 
	Clearly $\dm(\g_n)=3n/(d+1)+O(1)$ and by  Corollary~\ref{coro:equal},  $\mu(\g_n)=(1+o(1))\frac{4(d-3)\pi^2}{n^2}$.
	
	%For odd $d$, let $G_n$ be a $(d+2)$-regular graph of order $n$ and middle blocks $K_1+K_{d+1}+K_1$.
	%As $d$ is odd, it is possible 
	%to remove a $2$-factor from each block of $G_n$ and then obtain a $d$-regular graph $\g'_n$.
	%Then $\dm(\g'_n)=\dm(G_n)=3n/(d+3)+O(1)$, and
	%$\mu(\g'_n)\le\mu(G_n)=(1+o(1))\frac{(d+1)\pi^2}{n^2}$.

	For odd $d$, consider the graph $(K_1+K_{d+1}+K_1)_m$ with $m=\lfloor n/(d+3)\rfloor$. We remove a $2$-factor from each block in this graph and call the resulting graph $G_n$. As $d$ is odd, it is possible to modify the end blocks of $G_n$ to obtain a $d$-regular $n$-vertex graph $\g'_n$. Then $\dm(\g'_n)=3n/(d+3)+O(1)$.
	By Theorem~\ref{thm:H} and Corollary~\ref{coro:equal}, $\mu(\g'_n)=(1+o(1))\mu(G_n)\le(1+o(1))\frac{(d+1)\pi^2}{n^2}$.
	
	For even $d$, consider the graph $(K_1+K_{d-1}+K_2)_m$ with $m=\lfloor n/(d+2)\rfloor$. We remove a $1$-factor from each copy of $K_1+K_{d-1}+K_2$ and call the resulting graph $G_n$. Now we modify the end blocks of $G_n$ to obtain a $d$-regular $n$-vertex graph $\g'_n$. Then $\dm(\g'_n)=3n/(d+2)+O(1)$.  By Theorem~\ref{thm:H} and Corollary~\ref{coro:equal}, $\mu(\g'_n)=(1+o(1))\mu(G_n)\le(1+o(1))\frac{2(d-1)\pi^2}{n^2}$. 
	
	(iv) Consider the graph $(K_1+K_2+K_2)_m$ with $m=\lfloor n/5\rfloor$. We can modify the end blocks of this graph to obtain a graph $\g$ of order $n$ and $\delta=4$. Then $\dm(\g)=3n/5+O(1)$, and by Corollary~\ref{coro:equal}, $\mu(\g)=(1+o(1))\frac{4\pi^2}{n^2}$. 
This value is asymptotically smaller than the algebraic connectivity of graphs with $\delta=4$ obtained in virtue of Theorem~\ref{thm:muDiamter3nd+1}, which have diameter $3n/5+O(1)$ and $\mu=(1+o(1))\frac{3\pi^2}{n^2}$.
\end{proof}

We believe that the opposite direction of (b) should be true in general:
\begin{conjecture}\label{conj:asymp}
	For any $d\geq 3$ if $\Gamma_n$ is a sequence of graphs of $\delta=d$ (or a sequence of $d$-regular graphs) with asymptotically minimum  algebraic connectivity, then it has asymptotically maximum diameter that is
$(1+o(1))\frac{3n}{d+1}$.
\end{conjecture}

\section*{Acknowledgments}
 The first author was supported by a grant from IPM.
The second author carried out this work during a Humboldt Research Fellowship at the University of Hamburg. He thanks the Alexander von Humboldt-Stiftung for financial support.
The authors thank anonymous referees for several useful comments
which led to improvement of the paper's presentation.

\end{document}